\documentclass[2p]{article}
\usepackage{amssymb,amsmath,amsthm}
\usepackage{indentfirst}
\usepackage{exscale}
\usepackage{relsize}
\usepackage[numbers,sort&compress]{natbib}

\usepackage{geometry}
\usepackage{color}
\newcommand{\R}{\mathbb{R}}
\theoremstyle{definition}

\allowdisplaybreaks
\theoremstyle{remark}

\numberwithin{equation}{section}

\begin{document}

\title{\Large\bf{ Existence of ground state sign-changing solutions for a class of quasilinear scalar field equations originating from nonlinear optics }
 }
\date{}
\author {\ Xingyong Zhang$^{1,2}$\footnote{Corresponding author, E-mail address: zhangxingyong1@163.com}, {Xiaoli Yu}$^{1}$,  \\
{\footnotesize $^{1}$Faculty of Science, Kunming University of Science and Technology, Kunming, Yunnan, 650500, P.R. China.}\\
       {\footnotesize Kunming, Yunnan, 650500, P.R. China.}\\
       {\footnotesize $^{2}$Research Center for Mathematics and Interdisciplinary Sciences, Kunming University of Science and Technology,}\\
       {\footnotesize Kunming, Yunnan, 650500, P.R. China.}\\
 }
 \date{}
 \maketitle

 \begin{center}
 \begin{minipage}{15cm}

\small  {\bf Abstract:}
This paper is mainly concerned with the existence of ground state sign-changing solutions for a class of second order quasilinear elliptic equations in bounded domains which derived from nonlinear optics models. Combining a non-Nehari manifold method due to Tang and Cheng \cite{Tang and Cheng} and a quantitative deformation lemma with Miranda theorem, we obtain that the problem has at least one ground state sign-changing solution with two precise nodal domains. We also obtain that any weak solution of the problem has $C^{1,\sigma}$-regularity for some $\sigma\in(0,1)$. With the help of Maximum Principle, we reach the conclusion that  the energy of the ground state sign-changing solutions is strictly larger than twice that of the ground state solutions.

\end{minipage}
 \end{center}
  \allowdisplaybreaks
 \vskip2mm
\section{Introduction and main results }\label{section 1}
Self-trapped beams of light propagating in a nonlinear dielectric medium can be regard as special solutions of Maxwell's equations. In the case  that the magnetic field is everywhere transverse to the direction of propagation, the solutions are called TM-modes. Specifically, the study of the guided TM-modes is based on Maxwell's equations,
\begin{eqnarray}
\label{1}
\partial_t B=-c\nabla\wedge E,\;\;\;\;\;\;\partial_t D=-c\nabla\wedge H,
\end{eqnarray}
\begin{eqnarray}
\label{2}
\nabla \cdot B=0,\;\;\;\;\;\;\nabla \cdot D=0,
\end{eqnarray}
with the constitutive assumptions that
\begin{eqnarray}
\label{3}
 H=B ,\;\;\;\;\;\; D=\varepsilon(\langle E^2\rangle)E,
  \end{eqnarray}
  where $B$ is  the magnetic induction field,  $c$ is  the speed of light in a vacuum, $D$ is  the electric displacement field, $H$ is the magnetic field, $E$ is  the electric field,   $\langle E^2\rangle$ is the time-average of the intensity $E\cdot E$ of the electric field and $\varepsilon$  is a given function which satisfies the following conditions:\\
  $(1)$ $\varepsilon\in C([0,\infty),[0,\infty))\cap C^1([0,\infty),[0,\infty))$;\\
  $(2)$ $\varepsilon(0)>0$, $0<\varepsilon(\infty)=\lim_{s\rightarrow\infty}\varepsilon(s)<\infty $.\\
  More physical interest can be seen in \cite{Stuart12,Stuart13} and the standard procedure to study TM-modes in nonlinear optics  can be referred to \cite{Chen Y2,Chen Y3}. As is shown in  \cite{Stuart13}, the problem can be reduced to a single equation on $u$. In details, firstly, the function $\gamma:[0,\infty)\rightarrow[0,\infty)$ was defined by
 $
g(t)= \gamma(\dfrac{1}{2}t^2)t,
$
  where $g$ is the inverse of the function $f(s):=\varepsilon(\dfrac{1}{2}s^2)s$ with $s>0.$   If one assumes that $\gamma $ satisfies the following conditions:\\
  $(H_{1})$  $\gamma \in C([0,\infty),[0,\infty))\cap C^1([0,\infty),[0,\infty))$; \\
$(H_{2})$   $ \gamma'(t)\leq0<\gamma(t)+2t\gamma'(t);$\\
$(H_{3})$  $\gamma(0)>0$ and $\gamma(\infty)=\lim_{s\rightarrow\infty}\gamma(s)>0; $ \\
$(H_{4})$  $ \lim_{t\rightarrow0}t\gamma'(t)=0$;\\
$(H_{5})$ there exist constants $K>0$ and $ \sigma >0$ such that
$$
\lim_{t\rightarrow0}\dfrac{\gamma(t)-\gamma(0)}{t^\sigma}=-K.
$$
Then $(\ref{1})$,  $(\ref{2})$ and $(\ref{3})$ have a solution $u$ and $B$ has the form of $(\ref{3})$,
if and only if $u:[0,\infty)\rightarrow\mathbb{R}$  defined by
$$
u(r)=\dfrac{\omega}{kc}v(kr),
$$
 and $v:[0,\infty)\rightarrow\mathbb{R}$ satisfying
\begin{eqnarray}
\label{4}
\left(\gamma\left(\dfrac{1}{2}\left[v^2+\left(v'+\dfrac{v}{r}\right)^2\right]\right)\left(v'+\dfrac{v}{r}\right)\right)'
-\gamma\left(\dfrac{1}{2}\left[v^2+\left(v'+\dfrac{v}{r}\right)^2\right]\right)v+\lambda v=0,
\end{eqnarray}
where $(r,\theta,z)$ denotes cylindrical polar coordinate, $i_r$, $i_\theta  $, $i_z  $ denote curvilinear basis vectors, $B=u(r)\text{cos} (kz-\omega t)i_\theta$, $r\geq0$, $\lambda=\dfrac{\omega^2}{k^2c^2}$, $2\pi/k$ and $\omega$ denote the wavelength and frequency, respectively.
That is, finding the solutions of $(\ref{1})$, $(\ref{2})$ and $(\ref{3})$ is equivalent to finding the solutions of the second order differential equations (\ref{4}). Therefore, the magnetic field formulation can be described by a simple scalar second order differential equation. In \cite{Stuart C A},  the existence of nontrivial solutions for the nonlinear eigenvalue problem was established, which described self-trapped transverse magnetic field modes in a cylindrical optical fiber made from a self-focusing dielectric material. It is equal to establish the  existence of nontrivial solutions for a quasilinear second order differential equation in $H_0^1((0,\infty))$ and  the solutions are critical points of the  corresponding functional which has  mountain pass structure. In \cite{Zhou H S}, Stuart and Zhou  established  the existence of a minimizer of an integral functional, where the minimizer satisfies an Euler-Lagrange equation, and established the existence of guided TM-modes propagating through a self-focusing anisotropic dielectric.
\par
 Based on \cite{Stuart C A,Zhou H S}, Stuart \cite{C.A. Stuart} introduced the nonhomogeneous differential operator and established the result of the existence of at least two non-negative weak solutions for the following problem under the assumptions that $h\in L^2(\Omega),$ $h\geq 0$ and some other conditions,
\begin{eqnarray}
\label{af2}
 \begin{cases}
 -\text{div}\left(\phi\left(\dfrac{u^2+|\nabla u|^2}{2}\right)\nabla u\right)+\phi\left(\dfrac{u^2+|\nabla u|^2}{2}\right)u
=\lambda u+h,\;\;\;\;\;\; x\in \;\Omega ,\\
  u= 0,\;\;\;\;\;\; x\in \;\partial\Omega,
   \end{cases}
\end{eqnarray}
where $\Omega$ is a bounded domain of $\mathbb{R}^N,$ $N\geq 1$,  $\lambda \in \mathbb{R}$, $\Phi\in C^1([0,+\infty),\mathbb{R})$ and satisfes\\
$(\Phi_1)$  $\Phi(0)=0$, $\phi=\Phi'$ is non-increasing on $[0,+\infty)$ and $\phi(\infty)=\lim_{s\rightarrow\infty}\phi(s)>0;$ \\
$(\Phi_2)$ putting $\Phi(s^2)=g(s),$ there exists $\rho>0$ such that\\
$$ g(t)\geq g(s)+g'(s)(t-s)+\rho(t-s)^2,\;\;\;\;\mbox{for all}\; s,t\geq 0;$$
$(\Phi_3)$ $\lim_{s\rightarrow\infty} (\Phi(s)-\Phi'(s)s)<+\infty$.\\
$(\Phi2)$ is used to ensure the ellipticity of $(\ref{af2}).$ They also obtained that the equation $(\ref{af2})$ has a non-negative weak solution when  $h\equiv0$ and $\Phi$ satisfies $(\Phi_1)$ and $(\Phi_2)$. Subsequently, Jeanjean and R\u{a}dulescu \cite{Jeanjean L} extended the results in \cite{C.A. Stuart}. They studied the following second order quasilinear elliptic equation:
\begin{eqnarray}
\label{af3}
 \begin{cases}
 -\text{div}\left(\phi\left(\dfrac{u^2+|\nabla u|^2}{2}\right)\nabla u\right)+\phi\left(\dfrac{u^2+|\nabla u|^2}{2}\right)u
=f(u)+h,\;\;\;\;\;\; x\in \;\Omega ,\\
  u= 0,\;\;\;\;\;\; x\in \;\partial\Omega,
   \end{cases}
\end{eqnarray}
where $h\in L^2(\Omega)$ is non-negative, $\Omega$ is a bounded domain of $\mathbb{R}^N$ with $N\geq 1.$ They optimized the assumptions about the quasilinear part. Specifically,  they assumed that $\Phi$ satisfies the following conditions:\\
$(\phi_1)$ $\phi\in C([0,+\infty), \mathbb{R}^+)$ and there exist some constants $0<\phi_0\leq\phi_1$ such that
$$\phi_0 \leq\phi(s)\leq\phi_1,\;\;\;\;\;s\in[0,+\infty);$$
$(\phi_2)$ the function $s\mapsto \Phi(s^2)$ is convex on $\mathbb{R}.$ \\
They obtained that  problem $(\ref{af3})$ has a non-negative solution, when $h\gneqq 0$ and $f $ satisfies sublinear growth.  For $h\equiv0$,  they got the result that  the equation $(\ref{af3})$ has a  nontrivial solution and also obtain a nonexistence result. Meanwhile, when $ f$ satisfies linear growth, they showed the result that the problem $(\ref{af3})$  has at least one or two non-negative solutions.   The condition $(\phi_2)$ ensures the ellipticity of $(\ref{af3})$, and is weaker than $(\Phi_2) .$
In \cite{Pomponio A}, the following quasilinear elliptic problem was studied:
\begin{eqnarray}
\label{af4}
 \begin{cases}
 -\text{div}\left(\phi\left(\dfrac{u^2+|\nabla u|^2}{2}\right)\nabla u\right)+\phi\left(\dfrac{u^2+|\nabla u|^2}{2}\right)u
=g(u),\;\;\;\;\;\; x\in \;\mathbb{R}^N,\\
  u(x)\rightarrow 0,\;\;\;\;\;\; |x|\rightarrow +\infty,
   \end{cases}
\end{eqnarray}
where $N\geq 3.$ Combining the mountain pass theorem and some analytical skills, Pomponio and Watanabe got the existence of a non-negative  non-trivial  weak solution for (\ref{af4}),  if  the nonlinearity satisfies a variant of Berestycki-Lions' condition and $\phi$ satisfies $(\phi_1)$  and the following assumption: \\
$(\phi_2)'$ the function $s\mapsto \Phi(s^2)$ is strictly convex on $\mathbb{R}.$ \\
They got  a radial ground state solution $ u$ satisfying $u\in C^{1,\sigma}$ for some $ \sigma\in (0,1)$, if  the nonlinearity satisfies a variant of Berestycki-Lions' condition and $\phi$ satisfies $(\phi_1)$  and the following assumption: \\
$(\phi_3)$ $\phi\in C^1([0,+\infty), \mathbb{R}^+)$ and there exists $C>0$ such that
$$s|\phi'(s)|\leq C\;\;\;\mbox{and}\;\;\;\phi_0\leq\phi(s)+2s\phi'(s),\;\;\;\;\;\mbox{for all}\;s\in[0,+\infty).$$
They also proved that problem $(\ref{af4})$ has a ground state solution, when $(\phi_1)$, $(\phi_3)$ and the following condition holds:  \\
$(\phi_4)$ $\phi\in C^2([0,+\infty),\mathbb{R}^+)$ and there exists $C>0$ such that
$$s^2|\phi''(s)|\leq C,\;\;\;3\phi'(s)+2s|\phi''(s)|\leq 0\;\;\; \mbox{and} \;\;\;0\leq \phi(s)+5s\phi'(s)-2s^{2}|\phi''(s)|,\;\;\;\mbox{for all} \;s\in[0,+\infty).$$
\par
In the present paper, we investigate the existence of ground state sign-changing solutions for the following quasilinear elliptic problem:
\begin{eqnarray}
\label{af}
 \begin{cases}
 -\text{div}\left(\phi\left(\dfrac{u^2+|\nabla u|^2}{2}\right)\nabla u\right)+\phi\left(\dfrac{u^2+|\nabla u|^2}{2}\right)u
=f(u),\;\;\;\;\;\; x\in \;\Omega,\\
  u= 0,\;\;\;\;\;\; x\in \;\partial\Omega,
   \end{cases}
\end{eqnarray}
where $\Omega\subset\mathbb{R}^{N}$ is a bounded domain with a smooth boundary $\partial\Omega,$ $N=1, 2, 3$.  We assume $\Phi$ satisfies $(\phi_1)$, $(\phi_2)$ and the following condition:\\
$(\phi_5)$   $\phi\in C^1(\mathbb{R},  \mathbb{R}^+)$, $\phi$ is non-increasing, $|(\phi(s)s)'|\leq C$ for all $s\in[0,+\infty)$ and some $C>0$, and there exists a constant $\alpha\ge 0$ such that both $\phi\left(\frac{t^2}{2}\right)t^2$ and $\Phi\left(\frac{t^2}{2}\right)-\dfrac{1}{2(1+\alpha)}\phi\left(\frac{t^2}{2}\right)t^2$ are  convex on $t\in\mathbb{R}$.\\
Moreover, $f: \mathbb{R}\rightarrow\mathbb{R}$ satisfies the following assumptions:\\
$(F_{1})$ $f \in C (\mathbb{R},\mathbb{R})$ and there exist a $r>0$ and $\sigma_0\in[0,1)$ such that
$$
|f(t)|\le\sigma_0\phi_0(1+\lambda_1)|t|,\ \ \mbox{for all }|t|\le r;
$$
$(F_{2})$   there exist constants $C_{0}>0$ and $p\in (2,2^*)$ such that
$$
|f(t)|\leq C_{0}(1+|t|^{p-1}), \forall t\in\mathbb{R},
$$
where $2^*=6$ if $ N=3,$ and $2^*=+\infty$ if $N=1,2;$\\
$(F_{3})$ $\lim_{|t|\rightarrow \infty} \dfrac{f(t)}{t}=+\infty;$\\
 $(F_{4})$  there exists  $\theta_{0}\in (0,1)$ such that
$$
\left[\dfrac{f(\tau)}{\tau^{2\alpha+1}}-\dfrac {f(t\tau)}{(t\tau)^{2\alpha+1}}\right]\text{sign}(1-t)+\dfrac{\phi_0 (\theta_{0} \lambda_1+1 ) |1-t^{2\alpha}|}{(t\tau)^{2\alpha}}\geq 0, \mbox{ for all }  t>0 \mbox{ and all }\tau \in \mathbb{R}\backslash \{0\}
$$
or\\
$(F_{4})'$   there exists  $p_{0}> 0$ and $\gamma\geq p-2$ such that
$$
\left[\dfrac{f(\tau)}{\tau}-\dfrac {f(t\tau)}{t\tau}\right]\text{sign}(1-t)-p_0|\tau|^\gamma|1-t^\gamma|\geq 0, \mbox{ for all }  t>0 \mbox{ and all }\tau \in \mathbb{R}\backslash \{0\},
$$
where $\lambda_1$ is the first eigenvalue of $ (-\Delta, H_{0}^{1}(\Omega)),$ that is, $\lambda_1|u|^2\leq|\nabla u|^2$ for all $u\in H_{0}^{1}(\Omega)$, and $H^{1}_{0}(\Omega)$ denotes the usual Sobolev space equipped with the inner product and norm
$$
(u,v)=\int_{\Omega}\nabla u \nabla vdx,\;\;\;\|u\|=(\nabla u, \nabla u)^{1/2},\;\;\;\forall\;u, v\in H^{1}_{0}(\Omega).
$$
Throughout this paper, $L^{q}(\Omega)$ $(1\leq p<\infty)$ denotes the Lebesgue space with the norm $\|u\|_{q}=(\int_\Omega |u|^{q}dx)^{1/q}$.
It is well known that
 $ H^{1}_{0}(\Omega) \hookrightarrow L^{q}(\Omega)$ compactly with $ q\in [1,2^*)$, and  there exists a positive constant $C_q>0$ such that
\begin{eqnarray}
\label{b1}
\|u\|_q\leq C_q \|\nabla u\|_2.
\end{eqnarray}
\par
The corresponding functional of (\ref{af}) is defined on $H_{0}^{1}(\Omega)$ by
\begin{eqnarray}
\label{b2}
I(u)=\int_{\Omega}\Phi\left(\dfrac{u^2+|\nabla u|^2}{2}\right)dx-\int_{\Omega}F(u)dx,\;\;\;u\in H^{1}_{0}(\Omega),
\end{eqnarray}
where $F(s)=\int_0^sf(t)dt$.
\par
In recent years, the problem of sign-changing solutions for elliptic equations have been followed and many excellent results have been obtained. For example, the existence of sign-changing solutions for semilinear Schr\"{o}dinger equations was studied in \cite{Bartsch T1,Bartsch T2,Liu Z,Bartsch T3,Bartsch T4}. Meanwhile, the ground state sign-changing solutions of non-local Kirchhoff type problems also attracted some attentions. We refer the readers to \cite{ShuaiWei,Figueiredo G M,Tang and Cheng,Bartsch T,Cheng Chen} as some examples. Shuai \cite {ShuaiWei} studied the existence of  ground state sign-changing solutions for the following Kirchhoff-type problem by using a minimization argument and quantitative deformation lemma:
 \begin{eqnarray}
\label{ccc1}
 \begin{cases}
 -(a+b\int_\Omega|\nabla u|^2dx)\Delta u =f(u),\;\;\;\;\;\; x\in \;\Omega,\\
  u= 0,\;\;\;\;\;\; x\in \;\partial\Omega,
   \end{cases}
\end{eqnarray}
 where $a,b$ are positive constants, $\Omega$ is a bounded domain in $\mathbb R^N$ with $N=1,2,3$, and $f\in C^{1}$ satisfying the Nehari type monotonicity condition:\\
$(F)$  $\dfrac{f(t)}{t^{3}}$ is increasing on $(-\infty, 0) \cup (0,+\infty).$\\
 We must point out that the conditions $f\in C^{1}$ and $(F)$ play an important role in seeking ground state energy. For example, when the conditions $f\in C^{1}$ and $(F)$ are satisfied, He and Zou \cite{HeXZouW} also obtained the existence of ground state solutions for a class of  Kirchhoff problem by Nehari manifold's arguments, and Gasi\'{n}ski-Winkert \cite{Winkert P} and Liu-Dai  \cite{Liu W} obtained the existence of ground state sign-changing solutions for a class of double phase problem with nonlinear boundary condition. Under the conditions $f\in C^{1}$ and $\dfrac{f(t)}{t^{2m-1}}$ is increasing in $|t|>0,$  Figueiredo and Santos \cite{Figueiredo G M} extended the studies in  \cite {ShuaiWei} to  a generalized Kirchhoff equation.
\par
In order to weaken the assumptions on $f$,
Tang and Cheng \cite{Tang and Cheng} developed a direct non-Nehari manifold method and proved the existence of ground state sign-changing solutions of (\ref{ccc1}) under the weaker assumptions \\
$(F_1)'$:  $f(t)=o(t)$ as $ t\rightarrow 0$ and $f\in C^1$ actually in their proofs, \\
$(F_2)'$ there exist constants $\mathcal{C}_0>0$ and $p\in(4,2^*)$ such that \\
$$
|f(t) |\leq \mathcal{C}_0(1+|t|^{p-1}),\;\;\;\forall\;t\in \mathbb{R};
$$
$(F_3)'$ $ \lim_{|t|\rightarrow \infty}\dfrac{f(t)}{t^3}=\infty$;\\
$(F_{4})''$  there exists a $\theta_{0}\in (0,1)$ such that for all $t>0$ and $\tau \in \mathbb{R}$,
$$
\left[\dfrac{f(\tau)}{\tau^{3}}-\dfrac {f(t\tau)}{(t\tau)^{3}}\right]\text{sign}(1-t)+\dfrac{a \theta_{0} \lambda_1 |1-t^2|}{(t\tau)^2}\geq 0,
$$
where $\lambda_1$ is the first eigenvalue of $ (-\Delta, H_{0}^{1}(\Omega)).$ They also obtained that the energy of ground state sign-changing solutions is strictly larger than twice that of  ground state solutions. Subsequently, Cheng-Chen-Zhang \cite{Cheng Chen} combined the line in \cite{Tang and Cheng}  with a variant of  Miranda theorem to show that the existence of   ground state sign-changing solutions for a more general Kirchhoff-type Laplacian problem with the following assumption,\\
$(F_{4})'''$  there exists  $\theta_{0}\in (0,1)$  and $\alpha>0$ such that
$$
\left[\dfrac{f(\tau)}{\tau^{2\alpha+1}}-\dfrac {f(t\tau)}{(t\tau)^{2\alpha+1}}\right]\text{sign}(1-t)+\dfrac{a\theta_{0} \lambda_1 |1-t^{2\alpha}|}{(t\tau)^{2\alpha}}\geq 0, \mbox{ for all }  t>0 \mbox{ and all }\tau \in \mathbb{R}\backslash \{0\}.
$$
Especially, they  successfully weaken $f\in C^1 (\mathbb R,\mathbb R)$ to $f\in C(\mathbb R,\mathbb R)$.
\par
In this paper, following the line in \cite{Tang and Cheng,Cheng Chen}, we investigate the existence of ground state sign-changing solutions for the quasilinear elliptic problem $(\ref{af})$. The problem $(\ref{af})$ is different from the problems in  \cite{ShuaiWei,Figueiredo G M,Tang and Cheng,Bartsch T,Cheng Chen, Bartsch Liu Weth} because  both the quasilinear elliptic term and the nonlinear term depend on $u$ and $\nabla u$.
\par
 Let
 $$
 u^{+} :=\max \{u(x),0\}, \;\;\;\;\;u^{-} :=\min \{u(x),0\},
 $$
 $$
 \mathcal{M}:=\{u\in H_0^1(\Omega): u^\pm\neq0, \langle I'(u),u^{+}\rangle=\langle I'(u),u^{-}\rangle=0 \},
 $$
 and
 $$
\mathcal{N}:=\{u\in H_0^1(\Omega)\backslash\{0\},\langle I'(u),u\rangle=0\}.
$$
Obviously, $\mathcal{M}\subset \mathcal{N}.$  Next, we present our main results.
 \vskip2mm
 \noindent
 {\bf Theorem 1.1.}
If  $(\phi_1)$, $(\phi_2)$, $(\phi_5)$ with $\alpha>0$ and $(F_1)$-$(F_4)$ hold, then the problem $(\ref{af})$ possesses one ground state sign-changing solution $u_* \in \mathcal{M}$ such that $m_*:=I(u_*)=\inf_{u\in \mathcal{M}}I(u)$. Moreover,  $u_*$ has precisely two nodal domains and $u_* \in C^{1,\sigma}(\Omega)$ for some $\sigma\in(0,1)$.

\vskip2mm
 \noindent
 {\bf Theorem 1.2.}
If  $(\phi_1)$, $(\phi_2)$, $(\phi_5)$ with $\alpha>0$  and $(F_1)$-$(F_4)$ hold, then the problem $(\ref{af})$ has one ground state solution $u_0\in \mathcal{N}$ such that $u_0\in C^{1,\sigma}(\Omega)$ for some  $\sigma\in(0,1)$, and
\begin{eqnarray}\label{p1}
m_0:=I(u_0)=\inf_{\mathcal{N}}I\ge D_{\sigma_0}:=\dfrac{\alpha\phi_0(1-\theta_{0})}{2(1+\alpha)}\left(\dfrac{(1-\sigma_0)\phi_0r^{p-1}}{(C_0r^{p-1}+C_0)(C_p)^p}\right)^{2/(p-2)},
\end{eqnarray}
where $C_p$ is defined by (\ref{b1}).
\vskip2mm
 \noindent
 {\bf Theorem 1.3.} If  $(\phi_1)$, $(\phi_2)$, $(\phi_5)$ with $\alpha>0$ and $(F_1)$-$(F_4)$ hold, then $m_*> 2m_0$ and
 \begin{eqnarray}\label{p2}
\|u_*\|,\|u_0\| \ge \left(\dfrac{(1-\sigma_0)\phi_0r^{p-1}}{(C_0r^{p-1}+C_0)(C_p)^p}\right)^{1/(p-2)}.
\end{eqnarray}
\par
If $\alpha=0$ and we replace $(F_4)$ with $(F_4)'$, then we can also obtain the same conclusion as Theorem 1.1-Theorem 1.3. To be precise, we obtain the following result.
\vskip2mm
 \noindent
 {\bf Theorem 1.4.} If  $(\phi_1)$, $(\phi_2)$, $(\phi_5)$ with $\alpha=0$, $(F_1)$, $(F_2)$, $(F_3)$ and  $(F_4)'$ hold, then the problem $(\ref{af})$ possesses one ground state sign-changing solution $\hat{u}_*$ which has precisely two nodal domains, such that $\hat{m}_*:=I(\hat{u}_*)=\inf_{u\in \mathcal{M}}I(u)$, and a ground state solution $\hat{u}_0\in \mathcal{N}$ such that
\begin{eqnarray*}\label{p1}
\hat{m}_0&:=&I(u_0)=\inf_{\mathcal{N}}I\nonumber\\
&\geq& \hat{D}_{\sigma_0}:=\dfrac{p_0\gamma}{2(\gamma+2)}[\text{vol}(\Omega)]^{1-(\gamma+2)/p}
\left(\dfrac{(1-\sigma_0)\phi_0r^{p-1}}{C_0r^{p-1}+C_0}\right)^{(\gamma+2)/p}
\left(\dfrac{(1-\sigma_0)\phi_0r^{p-1}}{(C_0r^{p-1}+C_0)(C_p)^p}\right)^{[2(\gamma+2)]/[p(p-2)]},
\end{eqnarray*}
where  $\text{vol}(\Omega)=\int_\Omega 1dx<\infty$. Moreover,  $\hat{m}_*> 2\hat{m}_0$ and
 \begin{eqnarray}\label{p2}
\|\hat{u}\|,\|\hat{u}_0\| \ge \left(\dfrac{(1-\sigma_0)\phi_0r^{p-1}}{(C_0r^{p-1}+C_0)(C_p)^p}\right)^{1/(p-2)}.
\end{eqnarray}
 \par
\vskip2mm
 \noindent
 {\bf Remark 1.5.} In Theorem 1.1, we study the ground state sign-changing solutions which is not considered in \cite{Pomponio A}.
 Theorem 1.2 is different from that in \cite{Pomponio A} because we consider the case that $N=1,2,3$ rather than $N\ge 3$ and we only need that $\Phi(s)$ is convex rather than  strictly convex. What's more,  we only require that $\phi\in C^1(\mathbb{R}, \mathbb{R}^+)$ rather than $\phi\in C^2(\mathbb{R},\mathbb{R}^+)$ and $f $ satisfies the subcritical and super cubic growth rather than the Berestycki-Lions' type conditions. Different from \cite{ShuaiWei, Tang and Cheng, Cheng Chen}, the operator term in this paper is a function about $u$ and $\nabla u$, which makes the problem we studies more complex.

\vskip2mm
 \noindent
 {\bf Remark 1.6.} Obviously, $(F_1)$ is weaker than $(F_1)'$ even if $f\in C^1(\mathbb R,\mathbb R)$.  Moreover, in our results, by some careful calculation in the proofs, we  present the
  lower bounds of $\|u_*\|$ (or $\|\hat{u}_*\|$) and $\|u_0\| $  (or $\|\hat{u}_0\|$) and  a concrete lower bound  $D_{\sigma_0}$ (or $\hat{D}_{\sigma_0}$)  of the ground state energy  such that  $D_{\sigma_0}$ (or $\hat{D}_{\sigma_0}$)   is  non-decreasing on  $\sigma_0$, which are  not given in  \cite{ShuaiWei, Tang and Cheng, Cheng Chen}.
 \vskip2mm
 \noindent
  {\bf Remark 1.7.} There exist examples satisfying our assumptions. For example, let\\
$\bullet$ $\phi(s)\equiv 1$, for all $s\in\R$, which makes the problem $(\ref{af})$ become a second order semilinear elliptic equation;\\
$\bullet$  $\phi(s)=K+(1+s)^{-\beta}$ with $ \beta \in \left[0,\dfrac{2}{5}\right] $, $K\in (0,+\infty)$;\\
$\bullet$  $\phi(s)=K+(1+s)^{-\beta}$ with $ \beta \in \left[\dfrac{2}{5},\dfrac{1}{2}\right] $, $K\in\left [\dfrac{1}{2},+\infty\right)$.\\
Then it is easy to verify that  $(\phi_1), (\phi_2)$ and $(\phi_5)$ holds for all $\alpha\in [0,+\infty)$. Next, we give two examples of $f$ satisfying $(F_1)$-$(F_4)$. For example,  let \\
$(i)$ $f(t)=(|t|^{(2\alpha+1)}+\varrho|t|^{(2\alpha+1)/2})t,\;\forall\;t\in \mathbb{R}$, where $\alpha>0$ and $0<\varrho\leq2\sqrt{(1+\theta_0\lambda_1)\phi_0}$. (This example is similar to \cite[Appendix  A]{Cheng Chen});\\
$(ii)$ $f(t)=\beta t|t|^{2\alpha},\;  \forall\;t\in \mathbb{R}$, where $\alpha>0$ and $\beta >0$.\\
It is easy to see that the  example $(ii)$ satisfies $(F_1)$-$(F_4)$ but not satisfies the assumption $(A_3)$ in \cite{Cheng Chen}. The  example $(i)$ and the  example $(ii)$ also satisfy $(F_1)$-$(F_3)$ and $(F_4)'$.   Unfortunately, we do not know whether $(F_4)$ and  $(F_4)'$ are different.

\section{Preliminaries}\label{section 2}
The functional $I$ defined by (\ref{b2}) is well-defined on $H_{0}^{1}(\Omega)$. From $(\phi_1)$, $(\phi_2)$ and $(F_2)$, it holds that $I\in C^1(H_{0}^{1}(\Omega), \mathbb{R})$, and for any $u\in H_{0}^{1}(\Omega)$ and $\varphi\in H_{0}^{1}(\Omega),$
\begin{eqnarray}
\label{b3}
\langle I'(u),\varphi\rangle=\int_{\Omega}\phi\left(\dfrac{u^2+|\nabla u|^2}{2}\right)(u\varphi+\nabla u\cdot \nabla \varphi)dx
-\int_{\Omega}f(u)\varphi dx.
\end{eqnarray}
Let
$$
\Omega_u^{+} :=\{x\in\Omega|u(x)\geq0\}, \;\;\;\;\;\Omega_u^{-} :=\{x\in\Omega|u(x)<0\}.
$$
Then, for any $u\in H_0^1(\Omega)$, we have
\begin{eqnarray}
\label{b4}
\langle I'(u),u^{+}\rangle&=&\int_{\Omega}\phi\left(\dfrac{u^2+|\nabla u|^2}{2}\right)(uu^{+}+\nabla u\cdot \nabla u^{+})dx
-\int_{\Omega}f(u)u^{+} dx\nonumber\\
&=&\int_{\Omega_u^+\cup\Omega_u^-}\phi\left(\dfrac{u^2+|\nabla u|^2}{2}\right)(uu^{+}+\nabla u\cdot \nabla u^{+})dx
-\int_{\Omega_u^+\cup\Omega_u^-}f(u)u^{+} dx\nonumber\\
&=&\int_{\Omega_u^+}\phi\left(\dfrac{(u^+)^2+|\nabla u^+|^2}{2}\right)[(u^{+})^2+|\nabla u^{+}|^2]dx
-\int_{\Omega_u^+}f(u^+)u^{+} dx\nonumber\\
&&+\int_{\Omega_u^-}\phi\left(\dfrac{(u^-)^2+|\nabla u^-|^2}{2}\right)[u^{-}u^{+}+\nabla u^{-}\nabla u^{+}]dx
-\int_{\Omega_u^-}f(u^-)u^{+} dx\nonumber\\
&=&\int_{\Omega_u^+}\phi\left(\dfrac{(u^+)^2+|\nabla u^+|^2}{2}\right)[(u^{+})^2+|\nabla u^{+}|^2]dx
-\int_{\Omega_u^+}f(u^+)u^{+} dx\nonumber\\
&=&\int_{\Omega}\phi\left(\dfrac{(u^+)^2+|\nabla u^+|^2}{2}\right)[(u^{+})^2+|\nabla u^{+}|^2]dx
-\int_{\Omega}f(u^+)u^{+} dx\nonumber\\
&=&\langle I'(u^+),u^{+}\rangle.
\end{eqnarray}
Similarly,
\begin{eqnarray}
\label{b5}
\langle I'(u),u^{-}\rangle&=&\int_{\Omega}\phi\left(\dfrac{u^2+|\nabla u|^2}{2}\right)(uu^{-}+\nabla u\cdot \nabla u^{-})dx
-\int_{\Omega}f(u)u^{-} dx\nonumber\\
&=&\int_{\Omega_u^-}\phi\left(\dfrac{(u^-)^2+|\nabla u^-|^2}{2}\right)[(u^{-})^2+|\nabla u^{-}|^2]dx
-\int_{\Omega_u^-}f(u^-)u^{-} dx\nonumber\\
&=&\int_{\Omega}\phi\left(\dfrac{(u^-)^2+|\nabla u^-|^2}{2}\right)[(u^{-})^2+|\nabla u^{-}|^2]dx
-\int_{\Omega}f(u^-)u^{-} dx\nonumber\\
&=&\langle I'(u^-),u^{-}\rangle.
\end{eqnarray}
Consequently,
\begin{eqnarray}
\label{b6}
\langle I'(u),u\rangle=\langle I'(u^+),u^{+}\rangle+\langle I'(u^-),u^{-}\rangle.
\end{eqnarray}
Moreover, we also have
\begin{eqnarray}
\label{b7}
I(su^{+}+tu^{-})&=&\int_{\Omega}\Phi\left(\dfrac{(su^{+}+tu^{-})^2+|s\nabla u^{+}+t \nabla  u^{-}|^2}{2}\right)dx-\int_{\Omega}F(su^{+}+tu^{-})dx\nonumber\\
&=&\int_{\Omega_u^+\cup\Omega_u^-}\Phi\left(\dfrac{s^2(u^{+})^2+t^2(u^{-})^2+s^2|\nabla u^{+}|^2+t^2|\nabla  u^{-}|^2}{2}\right)dx-\int_{\Omega_u^+\cup\Omega_u^-}F(su^{+}+tu^{-})dx\nonumber\\
&=&\int_{\Omega_u^+}\Phi\left(\dfrac{s^2(u^{+})^2+s^2|\nabla u^{+}|^2}{2}\right)dx
+\int_{\Omega_u^-}\Phi\left(\dfrac{t^2(u^{-})^2+t^2|\nabla  u^{-}|^2}{2}\right)dx\nonumber\\
&&-\left[\int_{\Omega_u^+}F(su^{+})dx+\int_{\Omega_u^-}F(tu^{-})dx\right]\nonumber\\
&=&\int_{\Omega}\Phi\left(\dfrac{s^2(u^{+})^2+s^2|\nabla u^{+}|^2}{2}\right)dx
+\int_{\Omega}\Phi\left(\dfrac{t^2(u^{-})^2+t^2|\nabla  u^{-}|^2}{2}\right)dx\nonumber\\
&&-\left[\int_{\Omega}F(su^{+})dx+\int_{\Omega}F(tu^{-})dx\right]\nonumber\\
&=&I(su^{+})+I(tu^{-}),
\end{eqnarray}
where $s, t\geq 0.$\\

\vskip2mm
 \noindent
 \section{The regularity of weak solutions}\label{section 3}
 \vskip2mm
 \noindent
 {\bf Lemma 3.1.} {\it Assume that $(\phi_1)$, $(\phi_2)$, $(\phi_5)$ and $(F_1)$-$(F_4)$ hold. If $u$ is a weak solution of  (\ref{af})  then $u\in C^{1,\sigma}(\Omega)$ with some $\sigma\in(0,1)$.}\\
 {\bf Proof.} Following from \cite[Theorem 4.1]{Fan Zhao}, if the equations $(4.3)$-$(4.8)$ with $m(x)\equiv2$ and $r(x)\equiv p$ in \cite{Fan Zhao} hold, and $u$ is a weak solution, then $u\in L^\infty(\Omega)$. Next, we will claim every  weak solution of (\ref{af}) satisfies the equations $(4.3)$-$(4.8)$ in \cite{Fan Zhao}.
 \par
 Let
 $$E=\{u\in H^1_0(\Omega)\backslash \{0\}|I'(u)=0\},\;\;\;\; A(x,u,\nabla u)=\phi\left(\dfrac{u^2+|\nabla u|^2}{2}\right)\cdot\nabla u$$
  and
  $$B(x,u,\nabla u)=-\phi\left(\dfrac{u^2+|\nabla u|^2}{2}\right)u+f(u).$$
  Obviously, suppose that $u\in  E$. Then  $u$ is a weak solution of (\ref{af}).
  Following from $(\phi_1)$, it has
  $$
  A(x,u,\nabla u)\cdot\nabla u\geq\phi_0|\nabla u|^2,
  $$
  and
  $$
  |A(x,u,\nabla u)|\leq\phi_1|\nabla u|.
  $$
  From $(\phi_1)$, $(F_1)$ and $(F_2)$, it also has
  $$
   |B(x,u,\nabla u)|\leq\phi_1| u|+C_0| u|^{p-1}+C_0\leq C_1| u|^{p-1}+C_2
  $$
  where $C_1>0$ and $C_2>0$ are large constants. This shows that the equations $(4.3)$-$(4.8)$ in \cite{Fan Zhao}  hold in our paper. Then, we get that if $u\in E$,  $u\in L^\infty(\Omega).$ Similar to the discussion in part 4 of \cite{Pomponio A}, the equations $(4.3)$-$(4.8)$ in \cite{Pomponio A} also hold in our present paper.  Following from \cite[Proposition 4.1, Step 3]{Pomponio A}, we could use the regularity result of \cite[Chapter 4, Theorem 3.1, Theorem 5.2 and Theorem 6.2]{Ladyzhenskaya O A}  to obtain $u\in C^{1,\sigma}(\Omega)$ for some   $\sigma\in(0,1)$.\qed

\section{Proofs for ground state sign-changing solutions}\label{section 3}
\par
 \vskip2mm
 \noindent
 {\bf Lemma 4.1.} {\it Assume that $(\phi_1)$, $(\phi_2)$, $(\phi_5)$ and $(F_1)$-$(F_4)$ hold. Let $ u=u^{+}+u^{-}\in H_{0}^{1}(\Omega)$, $s,t\geq0$. Then }
 \begin{eqnarray}
\label{c1}
I(u)&\geq& I(su^{+}+tu^{-})+\dfrac{1-s^{2(1+\alpha)}}{2(1+\alpha)}\langle I'(u),u^{+}\rangle+\dfrac{1-t^{2(1+\alpha)}}{2(1+\alpha)}\langle I'(u),u^{-}\rangle\nonumber\\
&&+\dfrac{[\alpha-(\alpha+1)s^2+s^{2(1+\alpha)}]}{2(1+\alpha)}\phi_0(1-\theta_0)\|\nabla u^{+}\|_2^{2}
+\dfrac{{[\alpha-(\alpha+1)t^2+t^{2(1+\alpha)}]}}{{2(1+\alpha)}}\phi_0(1-\theta_0)\|\nabla u^{-}\|_2^{2}.
 \end{eqnarray}
 \par
\vskip2mm
 \noindent
 {\bf Proof.} It follows from $(F_4)$ that
 \begin{eqnarray}
\label{c2}
 &&\dfrac{1-t^{2(1+\alpha)}}{2(1+\alpha)}f(\upsilon)\upsilon+F(t\upsilon)-F(\upsilon)
 +\dfrac{(\theta_{0}\lambda_1+1)\phi_0}{2(1+\alpha)}[\alpha-(\alpha+1)t^2+t^{2(1+\alpha)}]\upsilon^{2}\nonumber\\
 &=&\int_{t}^{1}\left[\dfrac{f(\upsilon)}{\upsilon^{1+2\alpha}}-\dfrac{f(s\upsilon)}{(s\upsilon)^{1+2\alpha}}
 +\dfrac{\phi_0(\theta_{0}\lambda_1+1)(1-s^{2\alpha})}{(s\upsilon)^{2\alpha}}
\right]s^{1+2\alpha}\upsilon^{2(1+\alpha)}ds\geq0,\;\forall t\geq 0,\;\upsilon\in \mathbb{R}\backslash \{0\}.
 \end{eqnarray}
 Let $A(t)=\alpha-(\alpha+1)t^2+t^{2(1+\alpha)}$ and then $A(t)$ gets its minimum $A(1)=0$ and $1$ is unique minimum point, that is,
 \begin{eqnarray}
\label{hc3}
A(t)=\alpha-(\alpha+1)t^2+t^{2(1+\alpha)}\geq0,\;\;\mbox{for }\;t\geq0.
\end{eqnarray}
 By the fact that  $\phi$ is non-increasing and  the Mean value theorems for definite integrals, there exists a constant $\xi\in [a,b]$, where $[a,b]$ is an interval with $a\ge0$ and $b\ge0$, such that
 \begin{eqnarray}
 \label{c3}
\Phi(b)-\Phi(a)=\int_a^b\phi(s)ds=\phi(\xi)(b-a)\geq\phi(b)(b-a)
\end{eqnarray}
and
 \begin{eqnarray}
 \label{2c3}
\Phi(a)-\Phi(b)=\int_b^a\phi(s)ds=\phi(\xi)(a-b)\geq\phi(a)(a-b).
\end{eqnarray}
Thus if we set $b=\dfrac{u^2+|\nabla u|^2}{2}$ and $a=\dfrac{(su^{+}+tu^{-})^2+|s\nabla u^{+}+t \nabla  u^{-}|^2}{2}$,
 it follows from (\ref{c3}) that
 \begin{eqnarray}
 \label{3c3}
 &&\Phi\left(\dfrac{u^2+|\nabla u|^2}{2}\right)-\Phi\left(\dfrac{(su^{+}+tu^{-})^2+|s\nabla u^{+}+t \nabla  u^{-}|^2}{2}\right)\nonumber\\
 &\geq&\phi\left(\dfrac{u^2+|\nabla u|^2}{2}\right)\left[ \dfrac{u^2+|\nabla u|^2}{2}-\dfrac{(su^{+}+tu^{-})^2+|s\nabla u^{+}+t \nabla  u^{-}|^2}{2}\right].
\end{eqnarray}
If we set $b=\dfrac{(su^{+}+tu^{-})^2+|s\nabla u^{+}+t \nabla  u^{-}|^2}{2}$ and $a=\dfrac{u^2+|\nabla u|^2}{2}$,
it follows from (\ref{2c3}) that
 \begin{eqnarray}
 \label{4c3}
 &&\Phi\left(\dfrac{u^2+|\nabla u|^2}{2}\right)-\Phi\left(\dfrac{(su^{+}+tu^{-})^2+|s\nabla u^{+}+t \nabla  u^{-}|^2}{2}\right)\nonumber\\
 &\geq&\phi\left(\dfrac{u^2+|\nabla u|^2}{2}\right)\left[ \dfrac{u^2+|\nabla u|^2}{2}-\dfrac{(su^{+}+tu^{-})^2+|s\nabla u^{+}+t \nabla  u^{-}|^2}{2}\right].
\end{eqnarray}
 Then, using$(\ref{b2})$, $(\ref{b3})$-$(\ref{b7})$, $(\ref{c2})$, $(\ref{3c3})$, $(\ref{4c3})$,  $(\phi_1)$, $(\phi_2)$ and $(\phi_5)$, we obtain
\begin{eqnarray}
\label{c4}
&  & I(u)-I(su^{+}+tu^{-})\nonumber\\
& = &
\int_{\Omega}\left[\Phi\left(\dfrac{u^2+|\nabla u|^2}{2}\right)-\Phi\left(\dfrac{(su^{+}+tu^{-})^2+|s\nabla u^{+}+t \nabla  u^{-}|^2}{2}\right)\right]dx\nonumber\\
&&+\int_{\Omega}\left[F(su^{+}+tu^{-})-F(u)\right]dx\nonumber\\
&=&
\int_{\Omega}\int_{[(su^{+}+tu^{-})^2+|s\nabla u^{+}+t \nabla  u^{-}|^2]/2}^{[u^2+|\nabla u|^2]/2}\phi(\tau)d\tau dx
+\int_{\Omega}\left[F(su^{+}+tu^{-})-F(u)\right]dx\nonumber\\
&\ge &\int_{\Omega}\phi\left(\dfrac{u^2+|\nabla u|^2}{2}\right)\left[ \dfrac{u^2+|\nabla u|^2}{2}-\dfrac{(su^{+}+tu^{-})^2+|s\nabla u^{+}+t \nabla  u^{-}|^2}{2}\right]dx\nonumber\\
&&+\int_{\Omega}[F(su^{+})+F(tu^{-})-F(u^+)-F(u^-)]dx\nonumber\\
&\geq&\dfrac{1}{2}\int_{\Omega}\phi\left(\dfrac{u^2+|\nabla u|^2}{2}\right)\left(|u^{+}|^{2}+|u^{-}|^{2}+|\nabla u^{+}|^{2}+|\nabla u^{-}|^{2}-s^{2}|u^{+}|^{2}-t^{2}|u^{-}|^{2}-s^{2}|\nabla u^{+}|^{2}-t^{2}|\nabla u^{-}|^{2}\right)dx\nonumber\\
&&+\int_{\Omega}[F(su^{+})+F(tu^{-})-F(u^+)-F(u^-)]dx\nonumber\\
&=&\dfrac{1-s^{2(1+\alpha)}}{2(1+\alpha)}\int_{\Omega}\phi\left(\dfrac{u^2+|\nabla u|^2}{2}\right)(|u^{+}|^{2}+|\nabla u^{+}|^{2})dx+
\dfrac{1-t^{2(1+\alpha)}}{2(1+\alpha)}\int_{\Omega}\phi\left(\dfrac{u^2+|\nabla u|^2}{2}\right)(|u^{-}|^{2}+|\nabla u^{-}|^{2})dx\nonumber\\
&&+\dfrac{[\alpha-(\alpha+1)s^2+s^{2(1+\alpha)}]}{2(1+\alpha)}\int_{\Omega}\phi\left(\dfrac{u^2+|\nabla u|^2}{2}\right)(|u^{+}|^{2}+|\nabla u^{+}|^{2})dx\nonumber\\
&&+\dfrac{[\alpha-(\alpha+1)t^2+t^{2(1+\alpha)}]}{2(1+\alpha)}\int_{\Omega}\phi\left(\dfrac{u^2+|\nabla u|^2}{2}\right)(|u^{-}|^{2}+|\nabla u^{-}|^{2})dx\nonumber\\
&&+\int_{\Omega}[F(su^{+})+F(tu^{-})-F(u^+)-F(u^-)]dx\nonumber\\
&=&\dfrac{1-s^{2(1+\alpha)}}{2(1+\alpha)}\langle I'(u),u^{+}\rangle+\dfrac{1-t^{2(1+\alpha)}}{2(1+\alpha)}\langle I'(u),u^{-}\rangle
+\dfrac{1-s^{2(1+\alpha)}}{2(1+\alpha)}\int_{\Omega}f(u^{+})u^{+}dx+\dfrac{1-t^{2(1+\alpha)}}{2(1+\alpha)}\int_{\Omega}f(u^{-})u^{-}dx\nonumber\\
&&+\dfrac{[\alpha-(\alpha+1)s^2+s^{2(1+\alpha)}]}{2(1+\alpha)}\int_{\Omega}\phi\left(\dfrac{u^2+|\nabla u|^2}{2}\right)(|u^{+}|^{2}+|\nabla u^{+}|^{2})dx\nonumber\\
&&+\dfrac{[\alpha-(\alpha+1)t^2+t^{2(1+\alpha)}]}{2(1+\alpha)}\int_{\Omega}\phi\left(\dfrac{u^2+|\nabla u|^2}{2}\right)(|u^{-}|^{2}+|\nabla u^{-}|^{2})dx\nonumber\\
&&+\int_{\Omega}[F(su^{+})+F(tu^{-})-F(u^+)-F(u^-)]dx\nonumber\\
&=&\dfrac{1-s^{2(1+\alpha)}}{2(1+\alpha)}\langle I'(u),u^{+}\rangle+\dfrac{1-t^{2(1+\alpha)}}{2(1+\alpha)}\langle I'(u),u^{-}\rangle\nonumber\\
&&+\dfrac{[\alpha-(\alpha+1)s^2+s^{2(1+\alpha)}]}{2(1+\alpha)}\int_{\Omega}\phi\left(\dfrac{u^2+|\nabla u|^2}{2}\right)(|u^{+}|^{2}+|\nabla u^{+}|^{2})dx\nonumber\\
&&+\dfrac{[\alpha-(\alpha+1)t^2+t^{2(1+\alpha)}]}{2(1+\alpha)}\int_{\Omega}\phi\left(\dfrac{u^2+|\nabla u|^2}{2}\right)(|u^{-}|^{2}+|\nabla u^{-}|^{2})dx\nonumber\\
&&+\int_{\Omega}\left[\dfrac{1-s^{2(1+\alpha)}}{2(1+\alpha)}f(u^{+})u^{+}+F(su^{+})-F(u^+)\right]dx
+\int_{\Omega}\left[\dfrac{1-t^{2(1+\alpha)}}{2(1+\alpha)}f(u^{-})u^{-}+F(tu^{-})-F(u^-)\right]dx\nonumber\\
&\geq&\dfrac{1-s^{2(1+\alpha)}}{2(1+\alpha)}\langle I'(u),u^{+}\rangle+\dfrac{1-t^{2(1+\alpha)}}{2(1+\alpha)}\langle I'(u),u^{-}\rangle\nonumber\\
&&+\dfrac{[\alpha-(\alpha+1)s^2+s^{2(1+\alpha)}]}{2(1+\alpha)}\phi_0\int_{\Omega}(|u^{+}|^{2}+|\nabla u^{+}|^{2})dx
+\dfrac{[\alpha-(\alpha+1)t^2+t^{2(1+\alpha)}]}{2(1+\alpha)}\phi_0\int_{\Omega}(|u^{-}|^{2}+|\nabla u^{-}|^{2})dx\nonumber\\
&&-\dfrac{[\alpha-(\alpha+1)s^2+s^{2(1+\alpha)}]}{2(1+\alpha)}\phi_0\theta_{0}\int_{\Omega}|\nabla u^{+}|^{2}dx-\dfrac{[\alpha-(\alpha+1)t^2+t^{2(1+\alpha)}]}{2(1+\alpha)}\phi_0\theta_{0}\int_{\Omega}|\nabla u^{-}|^{2}dx\nonumber\\
&&+\int_{\Omega}\left[\dfrac{1-s^{2(1+\alpha)}}{2(1+\alpha)}f(u^{+})u^{+}dx+F(su^{+})-F(u^+)+\dfrac{[\alpha-(\alpha+1)s^2+s^{2(1+\alpha)}]}{2(1+\alpha)}\phi_0\theta_{0}\lambda_1|u^{+}|^{2}\right]dx\nonumber\\
&&+\int_{\Omega}\left[\dfrac{1-t^{2(1+\alpha)}}{2(1+\alpha)}f(u^{-})u^{-}+F(tu^{-})-F(u^-)+\dfrac{[\alpha-(\alpha+1)t^2+t^{2(1+\alpha)}]}{2(1+\alpha)}\phi_0\theta_{0}\lambda_1|u^{-}|^{2}\right]dx\nonumber\\
&=&\dfrac{1-s^{2(1+\alpha)}}{2(1+\alpha)}\langle I'(u),u^{+}\rangle+\dfrac{1-t^{2(1+\alpha)}}{2(1+\alpha)}\langle I'(u),u^{-}\rangle\nonumber\\
&&+\dfrac{[\alpha-(\alpha+1)s^2+s^{2(1+\alpha)}]}{2(1+\alpha)}\phi_0(1-\theta_{0})\int_{\Omega}|\nabla u^{+}|^{2}dx+\dfrac{[\alpha-(\alpha+1)t^2+t^{2(1+\alpha)}]}{2(1+\alpha)}\phi_0(1-\theta_{0})\int_{\Omega}|\nabla u^{-}|^{2}dx\nonumber\\
&&+\int_{\Omega}\left[\dfrac{1-s^{2(1+\alpha)}}{2(1+\alpha)}f(u^{+})u^{+}+F(su^{+})-F(u^+)
+\dfrac{[\alpha-(\alpha+1)s^2+s^{2(1+\alpha)}]}{2(1+\alpha)}\phi_0(\theta_{0}\lambda_1+1)|u^{+}|^{2}
\right]dx\nonumber\\
&&+\int_{\Omega}\left[\dfrac{1-t^{2(1+\alpha)}}{2(1+\alpha)}f(u^{-})u^{-}+F(tu^{-})-F(u^-)
+\dfrac{[\alpha-(\alpha+1)t^2+t^{2(1+\alpha)}]}{2(1+\alpha)}\phi_0(\theta_{0}\lambda_1+1)|u^{-}|^{2}
\right]dx\nonumber\\
&\geq&\dfrac{1-s^{2(1+\alpha)}}{2(1+\alpha)}\langle I'(u),u^{+}\rangle+\dfrac{1-t^{2(1+\alpha)}}{2(1+\alpha)}\langle I'(u),u^{-}\rangle+\dfrac{[\alpha-(\alpha+1)s^2+s^{2(1+\alpha)}]}{2(1+\alpha)}\phi_0(1-\theta_0)\|\nabla u^{+}\|_2^{2}\nonumber\\
&&+\dfrac{[\alpha-(\alpha+1)t^2+t^{2(1+\alpha)}]}{2(1+\alpha)}\phi_0(1-\theta_0)\|\nabla u^{-}\|_2^{2},\;\;\forall s,t\geq0.
\end{eqnarray}
 This shows that $(\ref{c1})$ holds.\qed
\par
\vskip2mm
 \noindent
 {\bf Corollary 4.2.} {\it Assume that $(\phi_1)$, $(\phi_2)$, $(\phi_5)$ and $(F_1)$-$(F_4)$ hold. If   $ u\in \mathcal{M}$, then}
 \begin{eqnarray}
\label{c5}
 I(u)&\geq& I(su^{+}+tu^{-})+\dfrac{[\alpha-(\alpha+1)s^2+s^{2(1+\alpha)}]}{2(1+\alpha)}\phi_0(1-\theta_0)\|\nabla u^{+}\|_2^{2}\nonumber\\
&&+\dfrac{[\alpha-(\alpha+1)t^2+t^{2(1+\alpha)}]}{2(1+\alpha)}\phi_0(1-\theta_0)\|\nabla u^{-}\|_2^{2},\;\;\forall s,t\geq0.
 \end{eqnarray}
 {\bf Proof.}  The result is obtained by Lemma 4.1 and the definition of $\mathcal{M}$ directly. \qed
\par
\vskip2mm
 \noindent
 {\bf Corollary 4.3.} {\it Assume that $(\phi_1)$, $(\phi_2)$, $(\phi_5)$ and $(F_1)$-$(F_4)$ hold. If   $ u\in \mathcal{M}$, then}
 \begin{eqnarray}
\label{c6}
 I(u)=\max_{s,t\geq0} I(su^{+}+tu^{-}).
 \end{eqnarray}
 {\bf Proof.} The proof is completed easily by observing (\ref{c5}). \qed
 \par
\vskip2mm
 \noindent
 {\bf Lemma 4.4.}
 {\it Assume that $(\phi_1)$, $(\phi_2)$, $(\phi_5)$ and $(F_1)$-$(F_4)$ hold. If $u=u^{+}+u^{-}\in H_{0}^{1}(\Omega)$ with $u^{\pm}\neq 0,$ then there exists a unique pair $(s_{u},t_{u})$ of positive numbers such that $s_{u}u^{+}+t_{u}u^{-}\in \mathcal{M}.$}
\par
\vskip2mm
 \noindent
 {\bf Proof.}
 We will use a method similar to \cite{Alves C O} to prove the existence of $(s_{u},t_{u}).$
For every $u\in H_{0}^{1}(\Omega)$ and $u^{\pm}\neq 0,$  let
\begin{eqnarray}
\label{c7}
\chi(s,t)&:=&\langle I'(su^{+}+tu^{-}),su^{+}+tu^{-}\rangle\nonumber\\
&=&\langle I'(su^{+}+tu^{-}),su^{+}\rangle+\langle I'(su^{+}+tu^{-}),tu^{-}\rangle\nonumber\\
&:=&(\chi_1(s,t),\chi_2(s,t)).\nonumber
\end{eqnarray}
By a direct calculation, it holds that
\begin{eqnarray}
\label{c8}
\chi_1(s,t)&=&\langle I'(su^{+}+tu^{-}),su^{+}\rangle\nonumber\\
&=&\int_{\Omega}\phi\left(\dfrac{(su^{+}+tu^{-})^2+|s\nabla u^{+}+t \nabla  u^{-}|^2}{2}\right)[su^{+}(su^{+}+tu^{-})+s\nabla u^{+}\cdot(s\nabla u^{+}+t\nabla u^{-})]dx\nonumber\\
&&-\int_{\Omega}f(su^{+}+tu^{-})su^{+}dx\nonumber\\
&=&\int_{\Omega_u^+}\phi\left(\dfrac{(su^{+})^2+s^2|\nabla u^{+}|^2}{2}\right)[(su^{+})^2+s^2|\nabla u^{+}|^2]dx-\int_{\Omega_u^+}f(su^{+})su^{+}dx\nonumber\\
&=&\int_{\Omega}\phi\left(\dfrac{(su^{+})^2+s^2|\nabla u^{+}|^2}{2}\right)[(su^{+})^2+s^2|\nabla u^{+}|^2]dx-\int_{\Omega}f(su^{+})su^{+}dx.
 \end{eqnarray}
Similarly,
\begin{eqnarray*}
\label{c9}
\chi_2(s,t)&=&\langle I'(su^{+}+tu^{-}),tu^{-}\rangle\nonumber\\
&=&\int_{\Omega}\phi\left(\dfrac{(su^{+}+tu^{-})^2+|s\nabla u^{+}+t \nabla  u^{-}|^2}{2}\right)[tu^{-}(su^{+}+tu^{-})+t\nabla u^{-}\cdot(s\nabla u^{+}+t\nabla u^{-})]dx\nonumber\\
&&-\int_{\Omega}f(su^{+}+tu^{-})tu^{-}dx\nonumber\\
&=&\int_{\Omega}\phi\left(\dfrac{(tu^{-})^2+t^2 |\nabla  u^{-}|^2}{2}\right)[(tu^{-})^2+t^2|\nabla u^{-}|^2]dx-\int_{\Omega}f(tu^{-})tu^{-}dx.
 \end{eqnarray*}
 Put $s>0$ small. Using $(\phi_1)$, $(F_1)$ and $(\ref{c8})$, and noting that $\sigma_0\in[0,1)$, we have
 \begin{eqnarray}
\label{c10}
\chi_1(s,s)&=&\int_{\Omega}\phi\left(\dfrac{(su^{+})^2+s^2|\nabla u^{+}|^2}{2}\right)[(su^{+})^2+s^2|\nabla u^{+}|^2]dx-\int_{\Omega}f(su^{+})su^{+}dx\nonumber\\
&\geq&\int_{\Omega}\phi_0[(su^{+})^2+s^2|\nabla u^{+}|^2]dx-\int_{\Omega}f(su^{+})su^{+}dx\nonumber\\
&\geq&\int_{\Omega}\phi_0[(su^{+})^2+s^2|\nabla u^{+}|^2]dx-\sigma_0\phi_0(1+\lambda_1)\int_{\Omega}(su^{+})^2dx\nonumber\\
&\geq&\phi_0(1+\lambda_1)\int_{\Omega}(su^{+})^2dx-\sigma_0\phi_0(1+\lambda_1)\int_{\Omega}(su^{+})^2dx\nonumber\\
&>&0.
\end{eqnarray}
Put $t>0$ large. From $(\phi_1)$, $(F_3)$ and $(\ref{c8})$,  there exists a $G_u>\frac{\phi_1\int_\Omega[|u^+|^2+|\nabla u^+|^2]dx}{\int_\Omega|u^+|^2dx}$ such  that
\begin{eqnarray}
\label{c11}
\chi_1(t,t)&=&\int_{\Omega}\phi\left(\dfrac{(tu^{+})^2+t^2|\nabla u^{+}|^2}{2}\right)[(tu^{+})^2+t^2|\nabla u^{+}|^2]dx-\int_{\Omega}f(tu^{+})tu^{+}dx\nonumber\\
&\leq&\int_{\Omega}\phi_1[(tu^{+})^2+t^2|\nabla u^{+}|^2]dx-G_u\int_{\Omega}(tu^{+})^2dx\nonumber\\
&<&0.
\end{eqnarray}
 Similarly, we have $\chi_2(s,s)>0$ for $s>0$ small and $\chi_2(t,t)<0$ for $t>0$ large. Thus, there exists $0<l<L$ such that
 \begin{eqnarray}
\label{c12}
\chi_1(l,l)>0,\;\;\;\; \chi_2(l,l)>0,\;\;\;\; \chi_1(L,L)<0,\;\;\;\; \chi_2(L,L)<0.
 \end{eqnarray}
From $(\ref{c8}),$ $(\ref{c10})$, $(\ref{c11})$ and $(\ref{c12}),$ we obtain that
\begin{eqnarray}
\label{c13}
\chi_1(l,t)&=&\int_{\Omega}\phi\left(\dfrac{(lu^{+})^2+l^2|\nabla u^{+}|^2}{2}\right)[(lu^{+})^2+l^2|\nabla u^{+}|^2]dx-\int_{\Omega}f(lu^{+})lu^{+}dx>0,\;\; \forall  t\in[l,L],
 \end{eqnarray}
 and
 \begin{eqnarray}
\label{c14}
\chi_1(L,t)&=&\int_{\Omega}\phi\left(\dfrac{(Lu^{+})^2+L^2|\nabla u^{+}|^2}{2}\right)[(Lu^{+})^2+L^2|\nabla u^{+}|^2]dx-\int_{\Omega}f(Lu^{+})Lu^{+}dx<0,\;\; \forall  t\in[l,L],
 \end{eqnarray}
since $\chi_1(l,t)$ and $\chi_1(L,t)$  are independent of $t$.
 Using the same steps, we have
\begin{eqnarray}
\label{c15}
\chi_2(s,l)>0,\;\;\;\; \chi_2(s,L)<0,\;\;\;\; \forall s\in[l,L].
 \end{eqnarray}
With the aid of  Miranda's Theorem \cite{Miranda C}, these imply that there exists  $(s_u,t_u)$ with $l<s_u,t_u<L$ such that $\chi_1(s_u,t_u)=\chi_2(s_u,t_u)=0.$ Then, $s_uu^++t_uu^-\in \mathcal{M}.$
\par
Next, we claim the uniqueness of $(s_u,t_u)$. Suppose that there exist two pairs $(s_1,t_1)$ and $(s_2,t_2)$ such that $s_iu^++t_iu^-\in \mathcal{M},$ $i=1,2.$ Obviously, $s_i\not=0$ and $t_i\not=0$, $i=1,2$. From Corollary 4.2, it holds that
\begin{eqnarray}
\label{c16}
I(s_1u^{+}+t_1u^{-})&\geq& I(s_2u^{+}+t_2u^{-})+\dfrac{\left[\alpha-(\alpha+1)\left(\dfrac{s_2}{s_1}\right)^2
+\left(\dfrac{s_2}{s_1}\right)^{2(1+\alpha)}\right]}{2(1+\alpha)}\phi_0(1-\theta_0)s_1^2\|\nabla u^{+}\|_2^{2}\nonumber\\
&&+\dfrac{\left[\alpha-(\alpha+1)\left(\dfrac{t_2}{t_1}\right)^2
+\left(\dfrac{t_2}{t_1}\right)^{2(1+\alpha)}\right]}{2(1+\alpha)}\phi_0(1-\theta_0)t_1^2\|\nabla u^{+}\|_2^{2}
 \end{eqnarray}
 and
 \begin{eqnarray}
\label{c17}
I(s_2u^{+}+t_2u^{-})&\geq& I(s_1u^{+}+t_1u^{-})
+\dfrac{\left[\alpha-(\alpha+1)\left(\dfrac{s_1}{s_2}\right)^2
+\left(\dfrac{s_1}{s_2}\right)^{2(1+\alpha)}\right]}{2(1+\alpha)}\phi_0(1-\theta_0)s_2^2\|\nabla u^{+}\|_2^{2}\nonumber\\
&&+\dfrac{\left[\alpha-(\alpha+1)\left(\dfrac{t_1}{t_2}\right)^2
+\left(\dfrac{t_1}{t_2}\right)^{2(1+\alpha)}\right]}{2(1+\alpha)}\phi_0(1-\theta_0)t_2^2\|\nabla u^{+}\|_2^{2}.
 \end{eqnarray}
 Note that $1$ the unique minimum point of $A(t)$. Then (\ref{hc3}), $(\ref{c16})$ and $(\ref{c17})$   yield $(s_1,t_1)=(s_2,t_2).$\qed
 \par
\vskip2mm
 \noindent
 {\bf Lemma 4.5.}
 {\it Assume that $(\phi_1)$, $(\phi_2)$, $(\phi_5)$ and $(F_1)$-$(F_4)$ hold. Then}
 $$
 \inf_{u\in \mathcal{M}}I(u)=m_* =\inf _{u\in H_0^1(\Omega), u^\pm\neq0}\max_{s,t\geq0}I(su^{+}+tu^{-}).
 $$
\par
\vskip2mm
 \noindent
 {\bf Proof.} Similar to the idea in \cite{Tang and Cheng}. From Corollary 4.3, one has
 \begin{eqnarray}
\label{c18}
\inf _{u\in H_0^1(\Omega), u^\pm\neq0}\max_{s,t\geq0}I(su^{+}+tu^{-})\leq \inf _{u\in \mathcal{M}}\max_{s,t\geq0}I(su^{+}+tu^{-})=\inf _{u\in \mathcal{M}}I(u)=m_*.
 \end{eqnarray}
Meanwhile, for any $u\in H_0^1(\Omega)$ with $ u^\pm\neq0,$ with the help of  Lemma 4.4, it holds  that
\begin{eqnarray*}
\label{c19}
\max_{s,t\geq0}I(su^{+}+tu^{-})\geq I(s_mu^++t_mu^-)\geq\inf _{u\in \mathcal{M}}I(u)=m_*.
 \end{eqnarray*}
This shows that
\begin{eqnarray}
\label{c20}
\inf _{u\in H_0^1(\Omega), u^\pm\neq0}\max_{s,t\geq0}I(su^{+}+tu^{-})\geq\inf _{u\in \mathcal{M}}I(u)=m_*.
 \end{eqnarray}
Following  from $(\ref{c18})$ and $(\ref{c20}),$ we obtain the conclusion.\qed
\par
\vskip2mm
 \noindent
  {\bf Lemma 4.6.} {\it Assume that $(\phi_1)$, $(\phi_2)$, $(\phi_5)$ and $(F_1)$-$(F_4)$ hold. Then $m_*=\inf_{u\in \mathcal{M}}I(u)\ge \inf_{u\in \mathcal{N}}I(u)=m_0\ge D_{\sigma_0}>0. $}\\
 {\bf Proof.} Obviously, $m*\ge m_0$ since $\mathcal{M}\subset\mathcal{N}$. It follows from  $(F_1)$ and $(F_2)$ that
 \begin{eqnarray*}
\label{cc21}
\int_{\Omega}f(u)u dx\leq\int_{\Omega}|f(u)u| dx\leq\sigma_0\phi_0(1+\lambda_1)\|u\|_2^2+\left(C_0+\frac{C_0}{r^{p-1}}\right)\|u\|^p_p.
 \end{eqnarray*}
 Since $\langle I'(u),u\rangle=0,\; \forall u\in \mathcal{N}$, together with  $(\phi_1)$ and $(\ref{b1})$, we have
 \begin{eqnarray}
\label{c21}
\phi_0\|\nabla u\|_2^2+\phi_0\| u\|_2^2&\leq&\int_{\Omega}\phi\left(\dfrac{u^2+|\nabla u|^2}{2}\right)(u^2+|\nabla u|^2)dx\nonumber\\
&  = & \int_{\Omega}f(u)u dx\nonumber\\
&\leq&\sigma_0\phi_0(1+\lambda_1)\|u\|_2^2+\left(C_0+\frac{C_0}{r^{p-1}}\right)\|u\|^p_p\nonumber\\
&\leq&\sigma_0\phi_0\|u\|_2^2+\sigma_0\phi_0\|\nabla u\|_2^2 +\left(\dfrac{C_0r^{p-1}+C_0}{r^{p-1}}\right)(C_p)^p\|\nabla u\|^p_2.
 \end{eqnarray}
Using $(\ref{c21})$, it holds that
  \begin{eqnarray}
\label{c22}
\|\nabla u\|_2\geq\left(\dfrac{(1-\sigma_0)\phi_0r^{p-1}}{(C_0r^{p-1}+C_0)(C_p)^p}\right)^{1/(p-2)}>0.
\end{eqnarray}
Similarly,
 \begin{eqnarray}
\label{ad2}
      \phi_0\|\nabla u\|_2^2+\phi_0\| u\|_2^2
&\leq &\int_{\Omega}\phi\left(\dfrac{u^2+|\nabla u|^2}{2}\right)(u^2+|\nabla u|^2)dx\nonumber\\
&  =  & \int_{\Omega}f(u)u dx\nonumber\\
&\leq & \sigma_0\phi_0(1+\lambda_1)\|u\|_2^2+\left(C_0+\frac{C_0}{r^{p-1}}\right)\|u\|^p_p\nonumber\\
&\leq & \sigma_0\phi_0\|u\|_2^2+\sigma_0\phi_0\|\nabla u\|_2^2 +\left(\dfrac{C_0r^{p-1}+C_0}{r^{p-1}}\right)\|u\|^p_p
 \end{eqnarray}
 which implies that
  \begin{eqnarray*}
\label{}
\|u\|_p^p\geq\dfrac{(1-\sigma_0)\phi_0r^{p-1}}{C_0r^{p-1}+C_0}\|\nabla u\|_2^2>0.
\end{eqnarray*}
Together with (\ref{c22}), we obtain
 \begin{eqnarray}
\label{ad1}
\|u\|_p\geq\left(\dfrac{(1-\sigma_0)\phi_0r^{p-1}}{C_0r^{p-1}+C_0}\right)^{1/p}
\left(\dfrac{(1-\sigma_0)\phi_0r^{p-1}}{(C_0r^{p-1}+C_0)(C_p)^p}\right)^{2/[p(p-2)]}>0.
\end{eqnarray}
   Note that for each $x\in \Omega$, there exists $\xi_{1,u,x}\in \left(0, \dfrac{u^2+|\nabla u|^2}{2}\right)$ such that
 \begin{eqnarray}
\label{c23}
\Phi\left(\dfrac{u^2+|\nabla u|^2}{2}\right)=\phi(\xi_{u,x} )\left(\dfrac{u^2+|\nabla u|^2}{2}\right),
\end{eqnarray}
from the Mean value theorem for definite integrals.
 Using  the conclusion of  $(\ref{c2})$ with $t=0$, we  have
  \begin{eqnarray}
\label{c24}
\dfrac{1}{2(1+\alpha)}f(\upsilon)\upsilon-F(\upsilon)+\dfrac{(\theta_{0}\lambda_1+1)\alpha\phi_0}{2(1+\alpha)}\upsilon^{2}\ge 0,
\;\;\;\; \forall \upsilon\in\mathbb{R}.
 \end{eqnarray}
Thus,  from $(\phi_1)$, $(\phi_2)$, $(\phi_5)$,  $(\ref{b2}),$ $(\ref{b3}),$  $(\ref{c22})$, $(\ref{c23})$ and $(\ref{c24})$, we obtain that
 \begin{eqnarray}
\label{c25}
I(u)&=&I(u)-\dfrac{1}{2(1+\alpha)}\langle I'(u),u\rangle\nonumber\\
&=&\int_{\Omega}\Phi\left(\dfrac{u^2+|\nabla u|^2}{2}\right)dx-\int_{\Omega}F(u)dx+\dfrac{1}{2(1+\alpha)}\int_{\Omega}f(u)u dx\nonumber\\
&&-\dfrac{1}{2(1+\alpha)}\int_{\Omega}\phi\left(\dfrac{u^2+|\nabla u|^2}{2}\right)(u^2+|\nabla u|^2)dx\nonumber\\
&=&\int_{\Omega}\int^{(u^2+|\nabla u|^2)/2}_{0}\phi(\tau)d\tau dx
-\dfrac{1}{2(1+\alpha)}\int_{\Omega}\phi\left(\dfrac{u^2+|\nabla u|^2}{2}\right)(u^2+|\nabla u|^2)dx
+ \int_{\Omega}\left[\dfrac{1}{2(1+\alpha)}f(u)u - F(u)\right]dx\nonumber\\
&=&\int_{\Omega}\phi(\xi_{u,x})\dfrac{u^2+|\nabla u|^2}{2}dx
-\dfrac{1}{2(1+\alpha)}\int_{\Omega}\phi\left(\dfrac{u^2+|\nabla u|^2}{2}\right)(u^2+|\nabla u|^2)dx
+ \int_{\Omega}\left[\dfrac{1}{2(1+\alpha)}f(u)u - F(u)\right]dx\nonumber\\
&\geq&\int_{\Omega}\phi\left(\dfrac{u^2+|\nabla u|^2}{2}\right)\dfrac{u^2+|\nabla u|^2}{2}dx
-\dfrac{1}{2(1+\alpha)}\int_{\Omega}\phi\left(\dfrac{u^2+|\nabla u|^2}{2}\right)(u^2+|\nabla u|^2)dx\nonumber\\
&&+ \int_{\Omega}\left[\dfrac{1}{2(1+\alpha)}f(u)u - F(u)\right]dx\nonumber\\
&\geq&\phi_0\int_{\Omega}\left[\dfrac{1}{2}(u^2+|\nabla u|^2)-\dfrac{1}{2(1+\alpha)}(u^2+|\nabla u|^2)\right]dx+ \int_{\Omega}\left[\dfrac{1}{2(1+\alpha)}f(u)u -F(u)\right]dx\nonumber\\
&\geq&\dfrac{\alpha}{2(1+\alpha)}\phi_0\|u\|_2^2+\dfrac{\alpha}{2(1+\alpha)}\phi_0\|\nabla u\|_2^2-\dfrac{(\theta_{0}\lambda_1+1)\alpha\phi_0}{2(1+\alpha)}\|u\|_2^2\nonumber\\
&\geq&\dfrac{\alpha\phi_0(1-\theta_{0})}{2(1+\alpha)}\|\nabla u\|_2^2\nonumber\\
&\geq&\dfrac{\alpha\phi_0(1-\theta_{0})}{2(1+\alpha)}\left(\dfrac{(1-\sigma_0)\phi_0r^{p-1}}{(C_0r^{p-1}+C_0)(C_p)^p}\right)^{2/(p-2)}\nonumber\\
&=&D_{\sigma_0}> 0,\ \ \forall u\in \mathcal{N}.
 \end{eqnarray}
 This implies $m_0\ge D_{\sigma_0}>0$.\qed
  \par
 \vskip2mm
 \noindent
 {\bf Lemma 4.7.} {\it Assume that $(\phi_5)$ holds. Let
 $$\Psi(u):=\displaystyle{\int_{\Omega}\phi\left(\dfrac{u^2+|\nabla u |^2}{2}\right)[u^2+|\nabla u|^2]dx }$$
 and
$$\Gamma(u):=\displaystyle{\int_{\Omega}\Phi\left(\dfrac{u^2+|\nabla u |^2}{2}\right)dx-\dfrac{1}{2(1+\alpha)}\int_{\Omega}\phi\left(\dfrac{u^2+|\nabla u |^2}{2}\right)[u^2+|\nabla u|^2]dx. }$$
Then both $\Psi$ and $ \Gamma$ are weakly sequentially lower semicontinuous, that is
$$\Psi(u)\leq\liminf_{n\rightarrow \infty}\Psi(u_n),\;\;\;\;\;
\;\;\;\;\;\Gamma(u)\leq\liminf_{n\rightarrow \infty}\Gamma(u_n),$$}
if $u_n\rightharpoonup u $ in $H_0^1(\Omega)$.\\
{\bf Proof.} The argument is similar to \cite[Lemma 3.2 (ii)]{Jeanjean L}. Setting $y=(u,\nabla u)$ and $z=(v,\nabla v)$ for every $u$, $v\in H_0^1(\Omega)$, one has $y\cdot z=uv+\nabla u\cdot \nabla v.$  Let
$$p(s):=\displaystyle{2s^2\phi(s^2) }.$$
Following from $(\phi_5)$ and similar to the arguments of (\ref{c3}) and (\ref{2c3}), we obtain that
\begin{eqnarray}
\label{D1}
\Psi(u)-\Psi(v)&=&\int_{\Omega}\left[p\left(\dfrac{|y|}{\sqrt{2}}\right)-p\left(\dfrac{|z|}{\sqrt{2}}\right)\right]dx\nonumber\\
&\geq&\int_{\Omega}p'\left(\dfrac{|z|}{\sqrt{2}}\right)\dfrac{|y|-|z|}{\sqrt{2}}dx\nonumber\\
&=&\int_{\Omega}\left[\sqrt{2}\phi'\left(\dfrac{|z|^2}{2}\right)|z|^3
+2\sqrt{2}\phi\left(\dfrac{|z|^2}{2}\right)|z|\right]\dfrac{|y|-|z|}{\sqrt{2}}dx\nonumber\\
&=&\int_{\Omega}\left[\phi'\left(\dfrac{|z|^2}{2}\right)|z|^2
+2\phi\left(\dfrac{|z|^2}{2}\right)\right][|z|(|y|-|z|)]dx.
\end{eqnarray}
From the fact that $(\phi_5)$, we know $\Psi(u)$ is Gateaux differentiable, and satisfies
\begin{eqnarray}
\label{D2}
\Psi'(v)(u-v)&=&\int_{\Omega}\phi'\left(\dfrac{v^2+|\nabla v |^2}{2}\right)(v^2+|\nabla v |^2)[v(u-v)+\nabla v \cdot \nabla(u-v)]dx\nonumber\\
&&+\int_{\Omega}2\phi\left(\dfrac{v^2+|\nabla v |^2}{2}\right)[v(u-v)+\nabla v \cdot \nabla(u-v)]dx\nonumber\\
&=&\int_{\Omega}\left[\phi'\left(\dfrac{|z|^2}{2}\right)|z|^2
+2\phi\left(\dfrac{|z|^2}{2}\right)\right][v(u-v)+\nabla v \cdot \nabla(u-v)]dx\nonumber\\
&=&\int_{\Omega}\left[\phi'\left(\dfrac{|z|^2}{2}\right)|z|^2
+2\phi\left(\dfrac{|z|^2}{2}\right)\right]z\cdot(y-z)dx.
\end{eqnarray}
The fact that $\Phi\left(\dfrac{t^2}{2}\right)$ is convex shows that
\begin{eqnarray}
\label{DD2}
 \phi\left(\dfrac{t^2}{2}\right)+\phi'\left(\dfrac{t^2}{2}\right)t^2\geq0,\;\;\;\mbox{for}\;\forall t\in \mathbb{R}.
\end{eqnarray}
 Together with $(\ref{D1})$, $(\ref{D2})$ and $(\ref{DD2})$, we get
 \begin{eqnarray}
\label{D3}
 \Psi(u)-\Psi(v)-\Psi'(v)(u-v)&\geq&\int_{\Omega}\left[\phi'\left(\dfrac{|z|^2}{2}\right)|z|^2
+2\phi\left(\dfrac{|z|^2}{2}\right)\right][|z|(|y|-|z|)-z\cdot(y-z)]dx\nonumber\\
 &=&\int_{\Omega}\left[\phi'\left(\dfrac{|z|^2}{2}\right)|z|^2
+2\phi\left(\dfrac{|z|^2}{2}\right)\right](|z||y|-|z|^2-z\cdot y+z^2])dx\nonumber\\
&=&\int_{\Omega}\left[\phi'\left(\dfrac{|z|^2}{2}\right)|z|^2
+2\phi\left(\dfrac{|z|^2}{2}\right)\right](|z||y|-z\cdot y])dx\geq0.
 \end{eqnarray}
 Setting $u=u_n$ and $v=u$ in $(\ref{D3})$, we have
 $$
 \Psi(u)\leq\Psi(u_n)-\Psi'(u)(u_n-u).
 $$
 Since $u_n\rightharpoonup u$ in $H_0^1(\Omega)$, then we have $\Psi'(u)(u_n-u)\rightarrow 0$.
 Thus, we reach the conclusion
 \begin{eqnarray}
\label{C1}
 \Psi(u)\leq\liminf_{n\rightarrow \infty}\Psi(u_n).
 \end{eqnarray}
 Similarly, let $h(s):=\Phi(s^2)-\dfrac{1}{2(1+\alpha)}[2s^2\phi(s^2)]$.  From $(\phi_1)$, $(\phi_2)$, $\phi\in C^1$, $\phi'\leq0$ and the fact that $h(s)$ is  convex, a direct  calculation implies
  \begin{eqnarray*}
\label{D4}
 &    &  \Gamma(u)-\Gamma(v)-\Gamma'(v)(u-v)\nonumber\\
 &\geq &
        \int_{\Omega}\left\{\phi\left(\dfrac{|z|^2}{2}\right)-\dfrac{1}{2(1+\alpha)}\left[2\phi\left(\dfrac{|z|^2}{2}\right)+\phi'\left(\dfrac{|z|^2}{2}\right)|z|^2
       \right]\right\}[|z|(|y|-|z|)-z\cdot(y-z)]dx\nonumber\\
 & =  & \int_{\Omega}\left[\dfrac{\alpha}{1+\alpha}\phi\left(\dfrac{|z|^2}{2}\right)-\dfrac{1}{2(1+\alpha)}\phi'\left(\dfrac{|z|^2}{2}\right)|z|^2
        \right][|z|(|y|-|z|)-z\cdot(y-z)]dx\nonumber\\
&  =  & \int_{\Omega}\left[\dfrac{\alpha}{1+\alpha}\phi\left(\dfrac{|z|^2}{2}\right)-\dfrac{1}{2(1+\alpha)}\phi'\left(\dfrac{|z|^2}{2}\right)|z|^2
        \right](|z||y|-z\cdot y)dx\nonumber\\
& \geq & 0.
 \end{eqnarray*}
Then, we also reach the conclusion that
\begin{eqnarray*}
\label{C2}
\Gamma(u)\leq\liminf_{n\rightarrow \infty}\Gamma(u_n).
\end{eqnarray*}
The proof is finished.\qed
 \par
 \vskip2mm
 \noindent
 {\bf Lemma 4.8.} {\it Assume that $(\phi_1)$, $(\phi_2)$, $(\phi_5)$ and $(F_1)$-$(F_4)$ hold. Then, there exists $u_*\in \mathcal{M} $ such that  $\inf_{u\in \mathcal{M}}I(u)=I(u_*)=m_*. $ }\\
{\bf Proof.} From Lemma 4.6, there exists a minimizing sequence $\{u_n\}\subset \mathcal{M}$ and $I(u_n)\rightarrow m_*$ as $n\rightarrow\infty$. Then by $(\ref{c25})$,  we have
 $$
 m_* +1\geq I(u_n)-\dfrac{1}{2(1+\alpha)}\langle I'(u_n),u_n\rangle
 \geq\dfrac{\alpha\phi_0(1-\theta_{0})}{2(1+\alpha)}\|\nabla u_n\|_2^2,
 $$
 for  $n\in \mathbb{N}$ large enough. This shows that $\|\nabla u_n\|_2^2$ is bounded and $\{u_n\}$ is bounded in $H_0^1(\Omega)$.  Then, there exists a subsequence, still denoted by $\{u_n\}$, such that
 \begin{eqnarray}\label{ad3}
 u_n\rightharpoonup u_*,  \;\;\;u_n^+\rightharpoonup u_1,\;\;\; u_n^-\rightharpoonup u_2,\;\;\;\mbox{in}\; H_0^1(\Omega);
 \end{eqnarray}
 \begin{eqnarray}\label{aa1}
 u_n\rightarrow u_*,\;\;\;u_n^+\rightarrow u_1,\;\;\; u_n^-\rightarrow u_2,\;\;\;\mbox{in}\; L^q(\Omega),\;q\in[1,2^*);
 \end{eqnarray}
  \begin{eqnarray}\label{ab2}
  u_n(x)\rightarrow u_*(x),\;\;\;u_n^+(x)\rightarrow u_1(x),\;\;\; u_n^-(x)\rightarrow u_2(x),\;\;\;\mbox{a.e.}\; x\in \Omega.
  \end{eqnarray}
  for some $u_*, u_1, u_2 \in H_0^1(\Omega)$. Note that for all $u\in H_0^1(\Omega)$, the maps $ u\rightarrow u^+$ and $ u\rightarrow u^-$ are continuous from $L^q(\Omega)$ into $L^q(\Omega)$ \cite[Lemma 2.3]{Castro A}. So (\ref{ad3}) implies that $u_n^+\rightharpoonup u^+$ and $u_n^-\rightharpoonup u^-$ in $H_0^1(\Omega)$ and (\ref{aa1}) implies that $u_n^+\to u_*^+$ and  $u_n^-\to u_*^-$ in $L^q(\Omega)$, and then
 \begin{eqnarray}
  \label{Z}
  u_n^+(x)\to u_*^+(x)\;\; \mbox{ and}\;\;  u_n^-(x)\to u_*^-(x) \mbox{ for a.e.}\;\; x\in\Omega .
  \end{eqnarray}
 By (\ref{ab2}) and the uniqueness of limit, we get
   $u_*^+(x)= u_1(x)\geq 0 $ and   $u_*^-(x)= u_2(x)\leq 0$ for a.e. $x\in \Omega$.
 Next, we can prove that $u_*\in \mathcal{M}.$ Indeed, by the fact that the embedding $H_0^1(\Omega)\hookrightarrow L^q(\Omega)$ is compact when $q\in[1,2^*)$, it is easy to know
\begin{eqnarray}
\label{c26}
u_n^\pm\rightarrow u_*^\pm\;\;\;\mbox{as}\;\;\;n\rightarrow\infty \;\;\mbox{in} \;L^p(\Omega),
\end{eqnarray}
for $p\in (2,2^*) $.
 Thus, one has
 \begin{eqnarray}\label{cd1}
 \int_\Omega|u_n^\pm|^pdx\rightarrow\int_\Omega|u_*^\pm|^pdx\;\;\;\mbox{as}\;n\rightarrow +\infty.
 \end{eqnarray}
 Since $\{u_n\}\subset \mathcal{M}\subset\mathcal{N}$, then, by (\ref{b4}) and (\ref{b5}), $\langle I'(u_n^\pm),u_n^\pm\rangle=\langle I'(u_n),u_n^\pm\rangle=0$. Similar to (\ref{c21}), it holds that
\begin{eqnarray}
\label{c33}
\|\nabla u_n^\pm\|_2\geq\left(\dfrac{(1-\sigma_0)\phi_0r^{p-1}}{(C_0r^{p-1}+C_0)(C_p)^p}\right)^{1/(p-2)}>0.
\end{eqnarray}
Similar to (\ref{ad2}) and together with (\ref{cd1}), one also has
\begin{eqnarray}
\label{c35}
\|u_*^\pm\|_p \geq\left(\dfrac{(1-\sigma_0)\phi_0r^{p-1}}{C_0r^{p-1}+C_0}\right)^{1/p}
\left(\dfrac{(1-\sigma_0)\phi_0r^{p-1}}{(C_0r^{p-1}+C_0)(C_p)^p}\right)^{2/p(p-2)}>0.
 \end{eqnarray}
This shows that $ u_*^\pm\neq 0.$

Note that  the  minimizing sequence $\{u_n\}\subset \mathcal{M}$. Then we have
\begin{eqnarray}
\label{c281}
\langle I'(u_n),u_n^+\rangle=\int_{\Omega}\phi\left(\dfrac{(u_n^+)^2+|\nabla u_n^+|^2}{2}\right)[(u_n^+)^2+|\nabla u_n^+|^2]dx-\int_{\Omega}f(u_n^+)u_n^+ dx=0
\end{eqnarray}
and
\begin{eqnarray}
\label{c282}
\langle I'(u_n),u_n^-\rangle=\int_{\Omega}\phi\left(\dfrac{(u_n^-)^2+|\nabla u_n^-|^2}{2}\right)[(u_n^-)^2+|\nabla u_n^-|^2]dx-\int_{\Omega}f(u_n^-)u_n^- dx=0.
\end{eqnarray}
Together with  $(F_2)$, there exists a positive constant $C_1$ such that
$$
|f(u_n^\pm)u_n^\pm|\leq C_0(|u_n^\pm|+|u_n^\pm|^p)\le C_1(1+|u_n^\pm|^p)
$$
and
$$
|F(u_n^\pm)|\leq C_0|u_n^\pm|+\frac{C_0}{p}|u_n^\pm|^p\le C_1(1+|u_n^\pm|^p).
$$
Combining \cite[Theorem A.2]{Willem},  there holds
$$
f(u_n^\pm)u_n^\pm\rightarrow f(u_*^\pm)u_*^\pm\;\;\;\mbox{as}\;\;\;n\rightarrow\infty \;\;\mbox{in} \;L^1(\Omega)
$$
and
$$
F(u_n^\pm)\rightarrow F(u_*^\pm)\;\;\;\mbox{as}\;\;\;n\rightarrow\infty \;\;\mbox{in} \;L^1(\Omega).
$$
Thus by Lebesgue dominated convergence theorem,
\begin{eqnarray}
\label{c31}
\lim_{n\rightarrow \infty}\int_{\Omega}f(u_n^\pm)u_n^\pm dx=\int_{\Omega}f(u_*^\pm)u_*^\pm dx
\end{eqnarray}
and
\begin{eqnarray}
\label{d32}
\lim_{n\rightarrow \infty}\int_{\Omega}F(u_n^\pm)dx=\int_{\Omega}F(u_*^\pm) dx.
\end{eqnarray}
Let
\begin{eqnarray*}
\Psi(u)=\int_{\Omega}\phi\left(\dfrac{u^2+|\nabla u |^2}{2}\right)[u^2+|\nabla u|^2]dx
\end{eqnarray*}
and
\begin{eqnarray*}
\Gamma(u)=\int_{\Omega}\Phi\left(\dfrac{u^2+|\nabla u |^2}{2}\right)dx-\dfrac{1}{2(1+\alpha)}\int_{\Omega}\phi\left(\dfrac{u^2+|\nabla u |^2}{2}\right)[u^2+|\nabla u|^2]dx.
\end{eqnarray*}
Following from Lemma 4.7, $(\ref{c281})$ and $(\ref{c31})$, it holds that
\begin{eqnarray*}
&&\int_{\Omega}\phi\left(\dfrac{(u_*)^2+|\nabla u_*^+|^2}{2}\right)[(u_*^+)^2+|\nabla u_*^+|^2]dx\nonumber\\
&\leq&\liminf_{n\rightarrow\infty} \int_{\Omega}\phi\left(\dfrac{(u_n^+)^2+|\nabla u_n^+|^2}{2}\right)[(u_n^+)^2+|\nabla u_n^+|^2]dx\nonumber\\
&=&\liminf_{n\rightarrow\infty}\int_{\Omega}f(u_n^+)u_n^+ dx\nonumber\\
&=& \int_{\Omega}f(u_*^+)u_*^+ dx.
\end{eqnarray*}
That is
 \begin{eqnarray}
\label{c36}
\langle I'(u_*),u_*^+\rangle\leq0.
\end{eqnarray}
Similarly, by (\ref{c282}) and (\ref{d32}), we also have
\begin{eqnarray}
\label{c37}
\langle I'(u_*),u_*^-\rangle\leq0.
\end{eqnarray}
 From $(\ref{b2}),$ $(\ref{b3}),$ $(\ref{c1}),$ $(\ref{c24})$, $(\ref{c26})$,  $(\ref{c36})$, $(\ref{c37})$,   Lemma 4.1, Lemma 4.7 and Fatou's Lemma, we have
 \\
\begin{eqnarray}
\label{c39}
m_*&=&\lim_{n\rightarrow \infty}\left[I(u_n)-\dfrac{1}{2(1+\alpha)}\langle I'(u_n),u_n\rangle\right]\nonumber\\
&=&\lim_{n\rightarrow \infty}\left[\int_{\Omega}\Phi\left(\dfrac{u_n^2+|\nabla u_n|^2}{2}\right)dx-\int_{\Omega}F(u_n)dx\right.\nonumber\\
&&\left.-\dfrac{1}{2(1+\alpha)}\int_{\Omega}\phi\left(\dfrac{u_n^2+|\nabla u_n|^2}{2}\right)(u_n^2+|\nabla u_n|^2)dx+\dfrac{1}{2(1+\alpha)}\int_{\Omega}f(u_n)u_n dx\right]\nonumber\\
&=&\lim_{n\rightarrow \infty}\left[\int_{\Omega}\Phi\left(\dfrac{u_n^2+|\nabla u_n|^2}{2}\right)dx
-\dfrac{1}{2(1+\alpha)}\int_{\Omega}\phi\left(\dfrac{u_n^2+|\nabla u_n|^2}{2}\right)(u_n^2+|\nabla u_n|^2)dx\right.\nonumber\\
&&\left.+ \int_{\Omega}\left[\dfrac{1}{2(1+\alpha)}f(u_n)u_n -F(u_n)\right]dx\right]\nonumber\\
&\geq&\liminf_{n\rightarrow \infty}\int_{\Omega}\left[\Phi\left(\dfrac{u_n^2+|\nabla u_n|^2}{2}\right)-\dfrac{1}{2(1+\alpha)}\phi\left(\dfrac{u_n^2+|\nabla u_n|^2}{2}\right)(u_n^2+|\nabla u_n|^2)-\dfrac{(\theta_{0}\lambda_1+1)\alpha\phi_0}{2(1+\alpha)}|u_n|^2\right]dx\nonumber\\
&&+ \liminf_{n\rightarrow \infty}\int_{\Omega}\left[\dfrac{1}{2(1+\alpha)}f(u_n)u_n - F(u_n)+\dfrac{(\theta_{0}\lambda_1+1)\alpha\phi_0}{2(1+\alpha)}|u_n|^2\right]dx\nonumber\\
&\geq&\liminf_{n\rightarrow \infty}\int_{\Omega}\left[\Phi\left(\dfrac{u_n^2+|\nabla u_n|^2}{2}\right)-\dfrac{1}{2(1+\alpha)}\phi\left(\dfrac{u_n^2+|\nabla u_n|^2}{2}\right)(u_n^2+|\nabla u_n|^2)\right]dx\nonumber\\
&&+\liminf_{n\rightarrow \infty}\int_{\Omega}-\dfrac{(\theta_{0}\lambda_1+1)\alpha\phi_0}{2(1+\alpha)}|u_n|^2dx
+ \liminf_{n\rightarrow \infty}\int_{\Omega}\left[\dfrac{1}{2(1+\alpha)}f(u_n)u_n - F(u_n)+\dfrac{(\theta_{0}\lambda_1+1)\alpha\phi_0}{2(1+\alpha)}|u_n|^2\right]dx\nonumber\\
&\geq&\int_{\Omega}\left[\Phi\left(\dfrac{u_*^2+|\nabla u_*|^2}{2}\right)dx-\dfrac{1}{2(1+\alpha)}\phi\left(\dfrac{u_*^2+|\nabla u_*|^2}{2}\right)(u_*^2+|\nabla u_*|^2)-\dfrac{(\theta_{0}\lambda_1+1)\alpha\phi_0}{2(1+\alpha)}u_*^{2}\right]dx\nonumber\\
&&+ \int_{\Omega}\left[\dfrac{1}{2(1+\alpha)}f(u_*)u_* - F(u_*)+\dfrac{(\theta_{0}\lambda_1+1)\alpha\phi_0}{2(1+\alpha)}u_*^{2}\right]dx\nonumber\\
&=&\int_{\Omega}\Phi\left(\dfrac{u_*^2+|\nabla u_*|^2}{2}\right)dx
-\dfrac{1}{2(1+\alpha)}\int_{\Omega}\phi\left(\dfrac{u_*^2+|\nabla u_*|^2}{2}\right)(u_*^2+|\nabla u_*|^2)dx
+ \int_{\Omega}\left[\dfrac{1}{2(1+\alpha)}f(u_*)u_* -F(u_*)\right]dx\nonumber\\
&=&I(u_*)-\dfrac{1}{2(1+\alpha)}\langle I'(u_*),u_*\rangle\nonumber\\
&\geq&\sup_{s,t\geq 0}\left[I(su_*^{+}+tu_*^{-})+\dfrac{1-s^{2(1+\alpha)}}{2(1+\alpha)}\langle I'(u_*),u_*^{+}\rangle+\dfrac{1-t^{2(1+\alpha)}}{2(1+\alpha)}\langle I'(u_*),u_*^{-}\rangle\right]-\dfrac{1}{2(1+\alpha)}\langle I'(u_*),u_*\rangle\nonumber\\
&=&\sup_{s,t\geq 0}\left[I(su_*^{+}+tu_*^{-})-\dfrac{s^{2(1+\alpha)}}{2(1+\alpha)}\langle I'(u_*),u_*^{+}\rangle-\dfrac{t^{2(1+\alpha)}}{2(1+\alpha)}\langle I'(u_*),u_*^{-}\rangle\right]\nonumber\\
&\geq&\sup_{s,t\geq 0}I(su_*^{+}+tu_*^{-})\geq m_*.
\end{eqnarray}
This implies
\begin{eqnarray}\label{aa2}
\sup_{s,t\geq 0}\left[I(su_*^{+}+tu_*^{-})-\dfrac{s^{2(1+\alpha)}}{2(1+\alpha)}\langle I'(u_*),u_*^{+}\rangle-\dfrac{t^{2(1+\alpha)}}{2(1+\alpha)}\langle I'(u_*),u_*^{-}\rangle\right]=m_*
\end{eqnarray}
and
 \begin{eqnarray}\label{aa3}
 \sup_{s,t\geq 0}I(su_*^{+}+tu_*^{-})=m_*.
 \end{eqnarray}
Since
$$\langle I'(u_*),u_*^\pm\rangle\leq0,$$
we get that $\langle I'(u_*),u_*^\pm\rangle=0.$ Otherwise, a contradiction with (\ref{aa2}) and (\ref{aa3}) appears.
 This shows that $u_*\in \mathcal{M}$, and then $I(u_*)=m_*$ by  (\ref{c39}).

\qed
\par
Next, we we will show that $u_*$ is a critical point of $I,$ i.e. $I'(u_*)=0,$ by quantitative deformation lemma and the arguments  in  \cite{Cheng Chen}.

 \par
\vskip2mm
 \noindent
 {\bf Lemma 4.9.} {\it Assume that $(\phi_1)$, $(\phi_2)$, $(\phi_5)$, $(F_1)$-$(F_4)$ hold. If $u_*\in \mathcal{M}$ and $I(u_*)=m_*,$ then $u_*$ is a critical point of $I$.}
 \par
\vskip2mm
 \noindent
 {\bf Proof.}
 By contradiction, we suppose that $I'(u_*)\neq 0.$ Then there exist $\delta >0$ and $\alpha >0$ such that
\begin{eqnarray*}
\label{c40}
u\in H_0^1(\Omega),\;\|u-u_*\|\leq3\delta\;\mbox{and}\; \|I'(u)\| \geq \alpha.
\end{eqnarray*}
Since $u_*\in \mathcal{M}$, following from Corollary 4.2, one knows that
\begin{eqnarray}
\label{c41}
 I(su_*^{+}+tu_*^{-})&\leq& I(u_*)-\dfrac{[\alpha-(\alpha+1)s^2+s^{2(1+\alpha)}]}{2(1+\alpha)}\phi_0(1-\theta_0)\|\nabla u_*^{+}\|_2^{2}\nonumber\\
 &&-\dfrac{[\alpha-(\alpha+1)t^2+t^{2(1+\alpha)}]}{2(1+\alpha)}\phi_0(1-\theta_0)\|\nabla u_*^{-}\|_2^{2}\nonumber\\
 &=&m_*-\dfrac{[\alpha-(\alpha+1)s^2+s^{2(1+\alpha)}]}{2(1+\alpha)}\phi_0(1-\theta_0)\|\nabla u_*^{+}\|_2^{2}\nonumber\\
 &&-\dfrac{[\alpha-(\alpha+1)t^2+t^{2(1+\alpha)}]}{2(1+\alpha)}\phi_0(1-\theta_0)\|\nabla u_*^{-}\|_2^{2},\;\;\;\forall s, t\geq 0.
\end{eqnarray}
Similar to $(\ref{c13})$,  $(\ref{c14})$ and  $(\ref{c15})$, there exist $r_1\in(0,1)$ and $R_1\in(1,+\infty)$ such that
\begin{eqnarray}
\label{c42}
\langle I'(r_1u_*^{+}+tu_*^{-}),r_1u_*^{+}\rangle>0,\;\;\;\langle I'(R_1u_*^{+}+tu_*^{-}),R_1u_*^{+}\rangle<0,\;\;\;\forall\; t\in[r_1,R_1],
\end{eqnarray}
and
\begin{eqnarray}
\label{c43}
\langle I'(su_*^{+}+r_1u_*^{-}),r_1u_*^{-}\rangle>0,\;\;\;\langle I'(su_*^{+}+R_1u_*^{-}),R_1u_*^{-}\rangle<0,\;\;\;\forall\; s\in[r_1,R_1].
\end{eqnarray}
\par
Let $D=(r_1,R_1)\times(r_1,R_1)\subset \mathbb{R}^2$ and $g(s,t)=su_*^{+}+tu_*^{-}.$ It follows from (\ref{c41}) and the fact $u_*^\pm\neq 0$ that
\begin{eqnarray}
\label{c44}
c:=\max_{(s, t)\in \partial D}I(su_*^{+}+tu_*^{-})<m_*.
\end{eqnarray}
\par
For $\varepsilon:=\min \{(m_*-c)/3, 1, \alpha\delta/8\},$ $S:=B(u_*,\delta)$ and $S_{2\delta}:=B(u_*,2\delta)$,   Lemma 2.3 in \cite{Willem} yields a deformation $\eta\in C([0,1]\times H_0^1(\Omega),  H_0^1(\Omega))$ such that\\
\par
(a). $\eta(1,u)=u$ if $u\notin I^{-1}([m_*-2\varepsilon, m_*+2\varepsilon]) \cap S_{2\delta}$;\\
\par
(b). $\eta\left(1,I^{-1}((-\infty, m_*+\varepsilon])\cap S\right)\subset I^{-1}((-\infty, m_*-\varepsilon]);$\\
\par
(c). $I(\eta(1,u))\leq I(u),\;\forall u\in H_0^1(\Omega);$\\
\par
(d). $\eta(t,\cdot)$ is a homeomorphism of $ H_0^1(\Omega)$, $\forall t\in[0,1].$
\vskip2mm
 \noindent
 \par
  Following from Corollary 4.3, $I(su_*^{+}+tu_*^{-})\leq I(u_*)=m_*<m_*+\varepsilon$ for $s,t\geq 0.$
From (b), for $su_*^{+}+tu_*^{-}\in I^{-1}((-\infty, m_*+\varepsilon])$ with any $s,t\geq0$ and $|s-1|^2+|t-1|^2<\dfrac{\delta^2}{\|u_*\|^2}$,   one has
 \begin{eqnarray}
\label{c45}
I(\eta(1,su_*^{+}+tu_*^{-}))\leq m_*-\varepsilon.
 \end{eqnarray}
 \vskip2mm
 \noindent
 Let $B:=\{(s,t): |s-1|^2+|t-1|^2<\dfrac{\delta^2}{\|u_*\|^2}\}.$ Following (\ref{hc3}), it holds
  \begin{eqnarray}
\label{fc46}
\inf_{(s,t)\in\bar{D}\backslash B}[B(s)+B(t)]:=\beta_0>0.
\end{eqnarray}
On the other hand, by (c), (\ref{c22}), (\ref{c41}), (\ref{fc46}) and the fact $u_*^\pm \neq 0$, one has
 \begin{eqnarray}
\label{c46}
I(\eta(1,su_*^{+}+tu_*^{-}))&\leq& I(su_*^{+}+tu_*^{-})\nonumber\\
 &\leq&m_*-\dfrac{[\alpha-(\alpha+1)s^2+s^{2(1+\alpha)}]}{2(1+\alpha)}\phi_0(1-\theta_0)\|\nabla u_*^{+}\|_2^{2}\nonumber\\
 &&-\dfrac{[\alpha-(\alpha+1)t^2+t^{2(1+\alpha)}]}{2(1+\alpha)}\phi_0(1-\theta_0)\|\nabla u_*^{-}\|_2^{2}\nonumber\\
 &\leq&m_*-\dfrac{\phi_0(1-\theta_0)}{{2(1+\alpha)}}\min\{\|\nabla u_*^{+}\|_2^{2},\|\nabla u_*^{-}\|_2^{2}\}[B(s)+B(t)] \nonumber\\
 &\leq&m_*-\dfrac{\phi_0(1-\theta_0)\beta_0}{{2(1+\alpha)}}\min\{\|\nabla u_*^{+}\|_2^{2},\|\nabla u_*^{-}\|_2^{2}\}\nonumber\\
 &&\forall s, t\geq 0, |s-1|^2+|t-1|^2\geq\dfrac{\delta^2}{\|u_*\|^2}.
 \end{eqnarray}
 From (\ref{c44}), (\ref{c45}) and  (\ref{c46}), it holds that
 \begin{eqnarray}
\label{c47}
\max_{(s, t)\in \bar{D}}I(\eta(1,su_*^{+}+tu_*^{-}))<m_*.
\end{eqnarray}
\par
Next, we prove that $\eta(1, g(D))\cap \mathcal{M}\neq \emptyset,$ which implies a contradiction with  the definition of $m_*.$ Let us define $h(s,t)=\eta(1, g(s,t))$ and
 \begin{eqnarray*}
\label{}
\Psi_1(s,t):=\langle I'(g(s,t)),g^+(s,t)\rangle =\langle I'(su_*^{+}+tu_*^{-}),su_*^{+}\rangle,\ \ \forall \ s,t\ge 0,
\end{eqnarray*}
 \begin{eqnarray*}
\label{}
\Psi_2(s,t):=\langle I'(g(s,t)),g^-(s,t) \rangle=\langle I'(su_*^{+}+tu_*^{-}),tu_*^{-}\rangle,\ \ \forall \ s,t\ge 0,
\end{eqnarray*}
\begin{eqnarray*}
\label{}
\Psi(s,t):=\left( \langle I'(h(s,t)),h^+(s,t)\rangle, \langle I'(h(s,t)),h^-(s,t)\rangle\right)\in \R^2, \ \ \forall \ s,t\ge 0.
\end{eqnarray*}
   Following from $(a)$ and the definition of $\varepsilon$, it is easy to obtain that
 \begin{eqnarray}
\label{c48}
h(s,t)=g(s,t)=su_*^{+}+tu_*^{-}
\end{eqnarray}
when $(s,t)\in\partial D$.
\par
It follows from $(\ref{c42})$, $(\ref{c43})$ and $(\ref{c48})$  that
 \begin{eqnarray}
\label{c49}
\Psi_1(r_1,t)>0,\;\;\;\Psi_1(R_1,t)<0,\;\;\;\forall\; t\in[r_1,R_1],
\end{eqnarray}
and
\begin{eqnarray}
\label{c50}
\Psi_2(s,r_1)>0,\;\;\;\Psi_2(s,R_1)<0,\;\;\;\forall\; s\in[r_1,R_1].
\end{eqnarray}
 Consequently, combining $(\ref{c49})$ and $(\ref{c50})$ with a variant of Miranda Theorem in \cite[Lemma 2.4]{Li G}, there exists $(\hat{s}, \hat{t})\in D,$ such that
\begin{eqnarray}
\label{c51}
 \Psi(\hat{s}, \hat{t})=(0,0)\Longleftrightarrow  \langle I'(h(\hat{s}, \hat{t})),h^\pm(\hat{s}, \hat{t})\rangle=0,
\end{eqnarray}
which shows that $ h(\hat{s}, \hat{t})\in \mathcal{M}.$ Moreover, note that $h(\hat{s}, \hat{t})=\eta(1, g(\hat{s}, \hat{t}))$ and $(\hat{s}, \hat{t})\in D$. Hence, $\eta(1, g(D))\cap \mathcal{M}\neq \emptyset$. \qed

 \par
\vskip2mm
 \noindent
{\bf Proof of Theorem 1.1.}   Lemma 4.9 shows  that $u_*$ is a critical point of $I$ and then a sign-changing ground state solution of problem $(\ref{af})$. Furthermore, by Lemma 3.1, $u_*\in C^{1,\sigma}(\Omega)$ for some $\sigma\in (0,1)$. Next, we show that $u_*$ has exactly two nodal domains by contradiction. Let $ u_*=u_1+u_2+u_3
 ,$ where
 \begin{eqnarray*}
\label{d1}
&  & u_1\geq 0,\;\;u_2\leq 0, \;\;\Omega_1 \cap \Omega_2=\emptyset,\;\;u_1|_{\Omega\backslash(\Omega_1\cup\Omega_2)}
 =u_2|_{\Omega\backslash(\Omega_1\cup\Omega_2)}
 =u_3|_{\Omega_1\cup\Omega_2}=0,\\
\label{d2}
&  & \Omega_1=\{x\in \Omega:u_1(x)>0\},\;\;\;\;\Omega_2=\{x\in \Omega:u_2(x)<0\},\;\;\;\;\Omega_3=\Omega\setminus(\Omega_1\cup\Omega_2),
  \end{eqnarray*}
   and $\Omega_1$, $\Omega_2$ are connected open subsets of $\Omega.$
Let $w=u_1+u_2$. Then $w^+=u_1$,  $w^-=u_2$ and $w^\pm\neq 0.$ Then, by Lemma 4.4, there exists a unique pair $(s_w, t_w)$ of positive numbers such that
\begin{eqnarray*}
\label{}
 s_ww^++t_ww^- \in \mathcal{M},
  \end{eqnarray*}
  equivalently,
  \begin{eqnarray*}
\label{}
 s_wu_1+t_wu_2 \in \mathcal{M}.
  \end{eqnarray*}
  Then
   \begin{eqnarray*}
\label{}
 I(s_wu_1+t_wu_2) \geq m_*.
  \end{eqnarray*}
  By the fact that $I'(u_*)=0$, we have
  \begin{eqnarray}
\label{N1}
\langle I'(u_*),w^+\rangle=\int_{\Omega}\phi\left(\dfrac{u_*^2+|\nabla u_*|^2}{2}\right)(u_*w^{+}+\nabla u_*\cdot\nabla w^{+})dx
-\int_{\Omega}f(u_*)w^{+} dx=0,
\end{eqnarray}
$$
\langle I'(u_*),w^-\rangle=\int_{\Omega}\phi\left(\dfrac{u_*^2+|\nabla u_*|^2}{2}\right)(u_*w^{-}+\nabla u_*\cdot\nabla w^{-})dx
-\int_{\Omega}f(u_*)w^{-} dx=0.
$$
Combining with (\ref{b4}), (\ref{b5}), (\ref{b6}) and (\ref{N1}), it holds that
  \begin{eqnarray}
\label{N2}
\langle I'(w),w^+\rangle&=&\int_{\Omega}\phi\left(\dfrac{(u_*-u_3)^2+|\nabla u_*-\nabla u_3|^2}{2}\right)[(u_*-u_3)w^{+}+(\nabla u_*-\nabla u_3)\cdot\nabla w^{+}]dx-\int_{\Omega}f(u_*-u_3)w^{+} dx\nonumber\\
&=&\int\limits_{\Omega_1\cup \Omega_2\cup\Omega_3}\phi\left(\dfrac{(u_*-u_3)^2+|\nabla u_*-\nabla u_3|^2}{2}\right)[(u_*-u_3)w^{+}+(\nabla u_*-\nabla u_3)\cdot\nabla w^{+}]dx\nonumber\\
&&-\int\limits_{\Omega_1\cup \Omega_2\cup\Omega_3 }f(u_*-u_3)w^{+} dx\nonumber\\
&=&\int_{\Omega_1}\phi\left(\dfrac{(w^{+})^2+|\nabla w^{+}|^2}{2}\right)[(w^{+})^2+|\nabla w^{+}|^2]dx-\int_{\Omega_1}f(w^+)w^{+} dx\nonumber\\
&=&\int_{\Omega}\phi\left(\dfrac{u_*^2+|\nabla u_*|^2}{2}\right)[(w^{+})^2+|\nabla w^{+}|^2]dx-\int_{\Omega}f(u_*)w^{+} dx\nonumber\\
&=&\langle I'(u_*),w^+\rangle=0.
  \end{eqnarray}
Similarly,\\
  \begin{eqnarray}
\label{N3}
\langle I'(w),w^-\rangle&=&\int_{\Omega}\phi\left(\dfrac{u_*^2+|\nabla u_*|^2}{2}\right)[(w^{-})^2+|\nabla w^{-}|^2]dx-\int_{\Omega}f(u_*)w^{-} dx\nonumber\\
&=&\langle I'(u_*),w^-\rangle=0.
\end{eqnarray}
\par
\vskip2mm
 \noindent
Similar to the argument of  (\ref{b6}) and (\ref{b7}),  using (\ref{c25}), (\ref{N2}), (\ref{N3}),  Lemma 4.1 and Lemma 4.5,  we get
  \begin{eqnarray*}
\label{}
 m_*=I(u_*)&=&I(u_*)-\dfrac{1}{2(1+\alpha)}\langle I'(u_*), u_*\rangle\nonumber\\
 &=&I(w)+I(u_3)-\dfrac{1}{2(1+\alpha)}\langle I'(w), w\rangle-\dfrac{1}{2(1+\alpha)}\langle I'(u_3), u_3\rangle\nonumber\\
 &\geq&\sup_{s,t\geq0}\left[I(sw^++tw^-)+\dfrac{1-s^{2(1+\alpha)}}{2(1+\alpha)}\langle I'(w),w^{+}\rangle+\dfrac{1-t^{2(1+\alpha)}}{2(1+\alpha)}\langle I'(w),w^{-}\rangle\right]\nonumber\\
 &&+I(u_3)-\dfrac{1}{2(1+\alpha)}\langle I'(w), w\rangle-\dfrac{1}{2(1+\alpha)}\langle I'(u_3), u_3\rangle\nonumber\\
 &=&\sup_{s,t\geq0}\left[I(sw^++tw^-)-\dfrac{s^{2(1+\alpha)}}{2(1+\alpha)}\langle I'(w),w^{+}\rangle-\dfrac{t^{2(1+\alpha)}}{2(1+\alpha)}\langle I'(w),w^{-}\rangle\right]+I(u_3)-\dfrac{1}{2(1+\alpha)}\langle I'(u_3), u_3\rangle\nonumber\\
 &=&\sup_{s,t\geq0}I(sw^++tw^-)+I(u_3)-\dfrac{1}{2(1+\alpha)}\langle I'(u_3), u_3\rangle\nonumber\\
 &\geq&\sup_{s,t\geq0}I(sw^++tw^-)+\dfrac{\alpha\phi_0(1-\theta_{0})}{2(1+\alpha)}\|\nabla u_3\|_2^2\nonumber\\
 &\geq& m_*+\dfrac{\alpha\phi_0(1-\theta_{0})}{2(1+\alpha)}\|\nabla u_3\|_2^2.
 \end{eqnarray*}
This shows that $u_3\equiv0$ and then $u_*$ has exactly two nodal domains.
\qed

\vskip2mm
 \noindent
 {\section{ Proofs for ground state solutions}}\label{section 5}
\par
\vskip2mm
 \noindent
 In this section, we will finish the proofs of Theorem 1.2. The process of proof is similar to Theorem 1.1. Firstly,  we can obtain Lemma 5.1-Lemma 5.6 below by the similar arguments of Lemma 4.1, Corollary 4.2, Corollary 4.3, Lemma 4.4,  Lemma 4.5 and Lemma 4.8, whose proofs are omitted. We only present  the proof of the fact that $u_0$ is a critical point of $I$, which is also essentially same as Lemma 4.9 with some modifications of details.
 \par
 \vskip2mm
 \noindent
 {\bf Lemma 5.1.} {\it Assume that $(\phi_1)$, $(\phi_2)$, $(\phi_5)$ and $(F_1)$-$(F_4)$ hold. For $\forall  u\in H_{0}^{1}(\Omega)$ and $\tau\geq0$, it holds that }
 \begin{eqnarray*}
\label{e1}
 I(u)&\geq& I(\tau u)+\dfrac{1-\tau^{2(1+\alpha)}}{2(1+\alpha)}\langle I'(u),u\rangle +\dfrac{[\alpha-(\alpha+1)\tau^2+\tau^{2(1+\alpha)}]}{2(1+\alpha)}\phi_0(1-\theta_0)\|\nabla u\|_2^{2}.
 \end{eqnarray*}
 \vskip2mm
 \noindent
 {\bf Corollary 5.2.} {\it Assume that $(\phi_1)$, $(\phi_2)$, $(\phi_5)$ and $(F_1)$-$(F_4)$ hold. For $\forall  u\in \mathcal{N}$, it holds that}
 \begin{eqnarray*}
\label{e2}
 I(u)&\geq& I(\tau u)+\dfrac{[\alpha-(\alpha+1)\tau^2+\tau^{2(1+\alpha)}]}{2(1+\alpha)}\phi_0(1-\theta_0)\|\nabla u\|_2^{2}.
 \end{eqnarray*}
 \par
\vskip2mm
 \noindent
 {\bf Corollary 5.3.} {\it Assume that $(\phi_1)$, $(\phi_2)$, $(\phi_5)$ and $(F_1)$-$(F_4)$ hold. If $ u\in \mathcal{N}$, then}
 \begin{eqnarray*}
\label{e3}
 I(u)=\max_{\tau\geq0} I(\tau u).
 \end{eqnarray*}
\vskip2mm
 \noindent
 {\bf Lemma 5.4.}
 {\it Assume that $(\phi_1)$, $(\phi_2)$, $(\phi_5)$ and $(F_1)$-$(F_4)$ hold. If $u\in H_{0}^{1}(\Omega)\backslash\{0\}$, then there exists a unique $\tau_{u}>0$ such that $\tau_{u}u\in \mathcal{N}.$}
  \par
\vskip2mm
 \noindent
 {\bf Lemma 5.5.}
 {\it Assume that $(\phi_1)$, $(\phi_2)$, $(\phi_5)$ and $(F_1)$-$(F_4)$ hold. Then
 $$
 \inf_{u\in \mathcal{N}}I(u)=:m_0 =\inf _{u\in H_0^1(\Omega), u\neq0}\max_{\tau\geq0}I(\tau u).
 $$}

\vskip2mm
 \noindent
 {\bf Lemma 5.6.} {\it  There exists a $u_0\in \mathcal{N}$ such that $I(u_0)=m_0\ge D_{\sigma_0}>0$.}

 \par
\vskip2mm
 \noindent
 {\bf Proof of Theorem 1.2.} We only need to prove that $u_0$ is  a critical point of $I$ and then  by Lemma 3.1, $u_0\in C^{1,\sigma}(\Omega)$ for some $\sigma\in (0,1)$. The proof is similar to Lemma 4.9 and we make the proof by contradiction.
 Assume that $I'(u_0)\neq 0$.  Then there exist constants  $\delta_0 >0$ and $\varrho_0 >0$ such that
\begin{eqnarray*}
\label{e4}
u\in H_0^1(\Omega),\;\|u-u_0\|\leq3\delta_0\;\Rightarrow \|I'(u)\|\geq \varrho_0.
\end{eqnarray*}
From Corollary 5.2, it holds that
\begin{eqnarray}
\label{e5}
 I(\tau u_0)&\leq& I(u_0)-\dfrac{[\alpha-(\alpha+1)\tau^2+\tau^{2(1+\alpha)}]}{2(1+\alpha)}\phi_0(1-\theta_0)\|\nabla u_0\|_2^{2}\nonumber\\
 &=&m_0-\dfrac{[\alpha-(\alpha+1)\tau^2+\tau^{2(1+\alpha)}]}{2(1+\alpha)}\phi_0(1-\theta_0)\|\nabla u_0\|_2^{2},\;\;\;\forall \tau\geq 0.
\end{eqnarray}
Similar to $(\ref{c10})$ and  $(\ref{c11})$, there exist $r_2\in(0,1)$ and $R_2\in(1,+\infty)$ such that
\begin{eqnarray}
\label{e6}
\langle I'(r_2u_0),r_2u_0\rangle>0,\;\;\;\langle I'(R_2u_0),R_2u_0\rangle<0.
\end{eqnarray}
Let $\Lambda=(r_2,R_2)\subset \mathbb{R}^+$ and $g(\tau )=\tau u_0.$ It follows from (\ref{e5}) that
\begin{eqnarray}
\label{e7}
c_0:=\max_{\tau=r_2 \;or\; R_2}I(\tau u_0)<m_0.
\end{eqnarray}
For $\varepsilon_0:=\min \{(m_0-c_0)/3, 1, \varrho_0\delta_0/8\},$ $S^0:=B(u_0,\delta_0),$ and $S_{2\delta_0}^0:=B(u_0,2\delta_0),$ \cite [Lemma 2.3]{Willem} yields a deformation $\eta_0\in C([0,1]\times H_0^1(\Omega), H_0^1(\Omega))$ such that\\
\par
$(a_0)$. $\eta_0(1,u)=u$ if $u\notin I^{-1}([m_0-2\varepsilon_0, m_0+2\varepsilon_0]) \cap S^0_{2\delta_0}$;\\
\par
$(b_0)$. $\eta_0\left(1,I^{-1}((-\infty, m_0+\varepsilon_0])\cap S^0\right)\subset I^{-1}((-\infty, m_0-\varepsilon_0]);$\\
\par
$(c_0)$. $I(\eta_0(1,u))\leq I(u),\;\forall u\in H_0^1(\Omega);$\\
\par
$(d_0)$. $\eta_0(t,\cdot)$ is a homeomorphism of $ H_0^1(\Omega)$, $\forall t\in[0,1].$
\vskip2mm
\par
  Following from Corollary 5.3, $I(\tau u_0)\leq I(u_0)=m_0<m_0+\varepsilon_0$ for $\tau \geq 0.$
From $(b_0)$, for $\tau u_0\in I^{-1}((-\infty, m_0+\varepsilon_0])$ with  $\tau\geq0$ and $|\tau -1|^2<\dfrac{\delta_0^2}{\|u_0\|^2}$,   one has
 \begin{eqnarray}
\label{e8}
I(\eta_0(1,\tau u_0))\leq m_0-\varepsilon_0.
 \end{eqnarray}
 Let $B_1:=\{\tau: |\tau -1|^2<\dfrac{\delta_0^2}{\|u_0\|^2}\}.$ Following (\ref{hc3}), it holds
  \begin{eqnarray}
\label{ffc46}
\inf_{\tau\in\bar{\Lambda}\backslash B_1}B(\tau):=\beta_1>0.
\end{eqnarray}
  Similar to (\ref{c21}),  (\ref{c22}), (\ref{ad2}) and (\ref{ad1}), one has
 $$
 \|\nabla u_0\|_2\geq\left(\dfrac{(1-\sigma_0)\phi_0r^{p-1}}{(C_0r^{p-1}+C_0)(C_p)^p}\right)^{1/(p-2)}>0
 $$
 and
 $$
 \| u_0\|_p\geq\left(\dfrac{(1-\sigma_0)\phi_0r^{p-1}}{C_0r^{p-1}+C_0}\right)^{1/p}
\left(\dfrac{(1-\sigma_0)\phi_0r^{p-1}}{(C_0r^{p-1}+C_0)(C_p)^p}\right)^{2/p(p-2)}>0.
 $$
On the other hand, by $(c_0)$, (\ref{e5}), (\ref{ffc46}) and the fact $u_0\neq0$, one has
 \begin{eqnarray}
\label{e9}
I(\eta_0(1,\tau u_0))&\leq& I(\tau u_0)\nonumber\\
 &\leq&m_0-\dfrac{[\alpha-(\alpha+1)\tau^2+\tau^{2(1+\alpha)}]}{2(1+\alpha)}\phi_0(1-\theta_0)\|\nabla u_0\|_2^{2}\nonumber\\
 &\leq&m_0-\dfrac{\phi_0(1-\theta_0)\beta_1}{2(1+\alpha)}\|\nabla u_0\|_2^{2},\;\;\;\forall \tau \geq 0, |\tau -1|^2\geq\dfrac{\delta_0^2}{\|u_0\|^2}.
 \end{eqnarray}
 From  (\ref{e7}), (\ref{e8}) and  (\ref{e9}), it holds that
 \begin{eqnarray}
\label{e10}
\max_{\tau\in \bar{\Lambda}}I(\eta_0(1,\tau u_0))<m_0.
\end{eqnarray}
\par
Next, we prove that $\eta_0(1, g(\Lambda))\cap \mathcal{N}\neq \emptyset,$ which implies a contradiction with the definition of $m_0.$ Let us define $h(\tau)=\eta_0(1, g(\tau))$ and
 \begin{eqnarray*}
\label{}
\psi(\tau) :=\langle I'(g(\tau)),g(\tau)\rangle=\langle I'(\tau u), \tau u\rangle,
\end{eqnarray*}
\begin{eqnarray*}
\label{}
\psi_0(\tau):=\langle I'(h(\tau)),h(\tau)\rangle.
\end{eqnarray*}
 It follows from $(a_0)$, (\ref{e7}) and the definition of $\varepsilon_0$ that
 \begin{eqnarray}
\label{e11}
h(\tau)=g(\tau)=\tau u_0,
\end{eqnarray}
when $\tau=r_2$ or $\tau=R_2$.
\par
In view of $(\ref{e6})$ and $(\ref{e11})$, we obtain that
 \begin{eqnarray}
\label{e12}
\psi_0(r_2)>0,\;\;\;\psi_0(R_2)<0.
\end{eqnarray}
 Consequently, combining $(\ref{e12})$  and  a variant of Miranda Theorem in \cite[Lemma 2.4]{Li G}, there exists $\tau_0\in \Lambda$ such that
\begin{eqnarray}
\label{e13}
 \psi_0(\tau_0)=0\Longleftrightarrow  \langle I'(h(\tau_0)),h(\tau_0)\rangle=0,
\end{eqnarray}
which shows that $h(\tau_0)\in \mathcal{M}$. Moreover, note that $h(\tau_0)=\eta(1, g(\tau_0)$ and $\tau_0\in \Lambda$. Hence, $\eta_0(1, g(\Lambda))\cap \mathcal{N}\neq \emptyset$.
\qed

\vskip2mm
\noindent
{\bf Proof of Theorem 1.3.}  Obviously, (\ref{p2}) holds by (\ref{c22}). Since
\begin{eqnarray}
\label{e14}
I(u_*)=m_*,\;\;\;\;\;I'(u_*)=0,\;\;\;\;\;u_*^{\pm}\neq0,
\end{eqnarray}
then it follows  from (\ref{b4}) and  (\ref{b5}) that $ u_*^+$, $ u_*^- \in \mathcal{N}$  and the definitions of $ \mathcal{N}$.
Combining with $(\ref{b2}),$ $(\ref{b7}),$ $(\ref{e14}),$ Corollary 4.3 and Lemma 5.5, we have
  \begin{eqnarray*}
\label{d5}
m_* & = & I(u_*) \nonumber\\
&  =  & I(u_*^{+}+u_*^{-})\nonumber\\
&  =  & I(u_*^{+})+I(u_*^{-})\nonumber\\
& \geq & 2 \inf_{u\in \mathcal{N}}I(u) \nonumber\\
&   =  & 2m_0.
\end{eqnarray*}
\par
 Suppose $ I(u_*)=m_*=2m_0$.  It holds that $ I (u_*^+)= I(u_*^-)=m_0$, since $I (u_*^+)\geq m_0$ and $ I(u_*^-) \geq m_0$. Following from Lemma 5.6 and the proof of Theorem 1.2, we get that  $ u_*^+$ is a critical point of $I$ and by Lemma 3.1, $u_*^+\in C^{1,\sigma}(\Omega)$. With the help of Maximum Principle \cite[Theorem 2.5.1]{Pucci P Serrin J B 2007}, we obtain $ u_*^+>0$ everywhere on $\Omega.$ This shows $ u_*^-=0$, which contraries to $u_*\in M.$ Thus, $ I(u_*)=m_*> 2m_0.$ \qed

 \vskip2mm
 \noindent
 {\section{ Proofs for  Theorem 1.4}}\label{section 5}
   We only present the proofs of Lemma 4.1-Lemma 4.4 and Lemma 4.6 with replacing $(F_4)$ with $(F_4)'$. The remainder proofs are similar to Theorem 1.1-Theorem 1.3.
  \par
\vskip2mm
 \noindent
{\bf Lemma 6.1.} {\it Assume that $(\phi_1)$, $(\phi_2)$, $(\phi_5)$,  $(F_1)$-$(F_3)$ and $(F_4)'$ hold. Let $ u=u^{+}+u^{-}\in H_{0}^{1}(\Omega)$, $s,t\geq0$ and $ \theta_0\in (0,1)$. Then }
 \begin{eqnarray}
\label{fc1}
 I(u)&\geq& I(su^{+}+tu^{-})+\dfrac{1-s^2}{2}\langle I'(u),u^{+}\rangle+\dfrac{1-t^2}{2}\langle I'(u),u^{-}\rangle\nonumber\\
 &&+\dfrac{p_0}{2(\gamma+2)}[\gamma-(\gamma+2)s^2+2s^{\gamma+2}]\| u^{+}\|_2^{\gamma+2}+\dfrac{p_0}{2(\gamma+2)}[\gamma-(\gamma+2)t^2+2t^{\gamma+2}]\| u^{-}\|_2^{\gamma+2}.
 \end{eqnarray}
 \par
\vskip2mm
 \noindent
 {\bf Proof.} It follows from $(F_4)'$ that
 \begin{eqnarray}
\label{fc2}
 &&\dfrac{1-t^2}{2}f(\upsilon)\upsilon+F(t\upsilon)-F(\upsilon)-\dfrac{p_0}{2(\gamma+2)}[\gamma-(\gamma+2)t^2+2t^{\gamma+2}]|\upsilon|^{\gamma+2}\nonumber\\
 &=&\int_{t}^{1}\left[\dfrac{f(\upsilon)}{\upsilon}-\dfrac{f(s\upsilon)}{s\upsilon}
 -p_0(1-s^\gamma)|\upsilon|^\gamma\right]s\upsilon^2ds\geq0,\;\;\forall t\geq 0,\;\;\;\upsilon\in \mathbb{R}\backslash \{0\}.
 \end{eqnarray}
  Let $B(t)=\gamma-(\gamma+2)t^2+2t^{\gamma+2}$. When  $t\in [0,1)$, $B'(t)\leq 0$,  and when $t>1$, $B'(t)> 0$. That is, $B(t)$ gets its minimum $B(1)=0$ and $1$ is the unique minimum point.
  Thus, it holds
 \begin{eqnarray}
\label{hhc3}
B(t)=\alpha-(\alpha+1)t^2+t^{2(1+\alpha)}\geq0,\;\;\mbox{for }\;t\geq0.
\end{eqnarray}
 Similar to (\ref{c4}), using $(\ref{b2})$, $(\ref{b3})$-$(\ref{b7})$, $(\ref{3c3})$, $(\ref{4c3})$, $(\ref{fc2})$,  $(\phi_1)$, $(\phi_2)$ and $(\phi_5)$, we get
 \begin{eqnarray}
\label{fc4}
&  & I(u)-I(su^{+}+tu^{-})\nonumber\\
& = &
\int_{\Omega}\left[\Phi\left(\dfrac{u^2+|\nabla u|^2}{2}\right)-\Phi\left(\dfrac{(su^{+}+tu^{-})^2+|s\nabla u^{+}+t \nabla  u^{-}|^2}{2}\right)\right]dx\nonumber\\
&&+\int_{\Omega}\left[F(su^{+}+tu^{-})-F(u)\right]dx\nonumber\\
&\geq&\dfrac{1}{2}\int_{\Omega}\phi\left(\dfrac{u^2+|\nabla u|^2}{2}\right)\left(|u^{+}|^{2}+|u^{-}|^{2}+|\nabla u^{+}|^{2}+|\nabla u^{-}|^{2}-s^{2}|u^{+}|^{2}-t^{2}|u^{-}|^{2}-s^{2}|\nabla u^{+}|^{2}-t^{2}|\nabla u^{-}|^{2}\right)dx\nonumber\\
&&+\int_{\Omega}[F(su^{+})+F(tu^{-})-F(u^+)-F(u^-)]dx\nonumber\\
&=&\dfrac{1-s^2}{2}\langle I'(u),u^{+}\rangle+\dfrac{1-t^2}{2}\langle I'(u),u^{-}\rangle\nonumber\\
&&+\dfrac{p_0}{2(\gamma+2)}[\gamma-(\gamma+2)s^2+2s^{\gamma+2}]\| u^{+}\|_2^{\gamma+2}
+\dfrac{p_0}{2(\gamma+2)}[\gamma-(\gamma+2)t^2+2t^{\gamma+2}]\| u^{-}\|_2^{\gamma+2}\nonumber\\
&&+\int_{\Omega}\left[\dfrac{1-s^2}{2}f(u^{+})u^{+}+F(su^{+})-F(u^+)-\dfrac{p_0}{2(\gamma+2)}[\gamma-(\gamma+2)s^2+2s^{\gamma+2}]| u^{+}|_2^{\gamma+2}\right]dx\nonumber\\
&&+\int_{\Omega}\left[\dfrac{1-t^2}{2}f(u^{-})u^{-}+F(tu^{-})-F(u^-)-\dfrac{p_0}{2(\gamma+2)}[\gamma-(\gamma+2)t^2+2t^{\gamma+2}]| u^{-}|_2^{\gamma+2}\right]dx\nonumber\\
&\geq&\dfrac{1-s^2}{2}\langle I'(u),u^{+}\rangle+\dfrac{1-t^2}{2}\langle I'(u),u^{-}\rangle\nonumber\\
&&+\dfrac{p_0}{2(\gamma+2)}[\gamma-(\gamma+2)s^2+2s^{\gamma+2}]\| u^{+}\|_2^{\gamma+2}
+\dfrac{p_0}{2(\gamma+2)}[\gamma-(\gamma+2)t^2+2t^{\gamma+2}]\| u^{-}\|_2^{\gamma+2},\;\;\forall s,t\geq0.\nonumber
 \end{eqnarray}
 This shows that (\ref{fc1}) holds.\qed

\vskip2mm
 \noindent
 {\bf Corollary 6.2.} {\it Assume that $(\phi_1)$, $(\phi_2)$, $(\phi_5)$, $(F_1)$-$(F_3)$  and $(F_4)'$ hold. If   $ u\in \mathcal{M}$, then}
 \begin{eqnarray*}
 I(u)&\geq& I(su^{+}+tu^{-})+\dfrac{p_0}{2(\gamma+2)}[\gamma-(\gamma+2)s^2+2s^{\gamma+2}]\| u^{+}\|_2^{\gamma+2}\nonumber\\
&&+\dfrac{p_0}{2(\gamma+2)}[\gamma-(\gamma+2)t^2+2t^{\gamma+2}]\| u^{-}\|_2^{\gamma+2}
 \end{eqnarray*}
for any $ s$, $t\geq0.$\\

  \par
\vskip2mm
 \noindent
 {\bf Corollary 6.3.} {\it Assume that $(\phi_1)$, $(\phi_2)$, $(\phi_5)$ $(F_1)$-$(F_3)$  and $(F_4)'$ hold. If   $ u\in \mathcal{M}$, then}
 \begin{eqnarray*}
 I(u)=\max_{s,t\geq0} I(su^{+}+tu^{-}).
 \end{eqnarray*}

 \par
\vskip2mm
 \noindent
 {\bf Lemma 6.4.}
 {\it Assume that $(\phi_1)$, $(\phi_2)$, $(\phi_5)$, $(F_1)$-$(F_3)$  and $(F_4)'$ hold. If $u=u^{+}+u^{-}\in H_{0}^{1}(\Omega)$ with $u^{\pm}\neq 0,$ then there exists a unique pair $(s_{u},t_{u})$ of positive numbers such that $s_{u}u^{+}+t_{u}u^{-}\in \mathcal{M}.$}
\par
\vskip2mm
 \noindent
 {\bf Proof.} The proof of the existence of $(s_{u},t_{u})$ is the same as that in Lemma 4.4, so we will omit it here.
Next, we show the uniqueness of $(s_u,t_u)$. Suppose that there exist two pairs $(s_1,t_1)$ and $(s_2,t_2)$ such that $s_iu^++t_iu^-\in \mathcal{M},$ $i=1,2.$ From Corollary 6.2, it holds that
\begin{eqnarray}
\label{fc16}
I(s_1u^{+}+t_1u^{-})&\geq& I(s_2u^{+}+t_2u^{-})+\dfrac{p_0}{2(\gamma+2)}\left[\gamma-(\gamma+2)\left(\dfrac{s_2}{s_1}\right)^2
+2\left(\dfrac{s_2}{s_1}\right)^{\gamma+2}\right]{s_1}^{\gamma+2}\|u^{+}\|_2^{\gamma+2}\nonumber\\
&&+\dfrac{p_0}{2(\gamma+2)}\left[\gamma-(\gamma+2)\left(\dfrac{t_2}{t_1}\right)^2
+2\left(\dfrac{t_2}{t_1}\right)^{\gamma+2}\right]{t_1}^{\gamma+2}\|u^{+}\|_2^{\gamma+2},
 \end{eqnarray}
 and
 \begin{eqnarray}
\label{fc17}
I(s_2u^{+}+t_2u^{-})&\geq& I(s_1u^{+}+t_1u^{-})+\dfrac{p_0}{2(\gamma+2)}\left[\gamma-(\gamma+2)\left(\dfrac{s_1}{s_2}\right)^2
+2\left(\dfrac{s_1}{s_2}\right)^{\gamma+2}\right]{s_2}^{\gamma+2}\|u^{+}\|_2^{\gamma+2}\nonumber\\
&&+\dfrac{p_0}{2(\gamma+2)}\left[\gamma-(\gamma+2)\left(\dfrac{t_1}{t_2}\right)^2
+2\left(\dfrac{t_1}{t_2}\right)^{\gamma+2}\right]{t_2}^{\gamma+2}\|u^{+}\|_2^{\gamma+2}.
 \end{eqnarray}
 Note that $1$ the unique minimum point of $B(t)$. Then (\ref{hhc3}), $(\ref{fc16})$ and $(\ref{fc17})$  yield $(s_1,t_1)=(s_2,t_2).$
 \qed
 \par
\vskip2mm
 \noindent
  {\bf Lemma 6.5.} {\it Assume that $(\phi_1)$, $(\phi_2)$, $(\phi_5)$, $(F_1)$-$(F_3)$  and $(F_4)'$ hold. Then $\hat{m}_*=\inf_{u\in \mathcal{M}}I(u)\ge \inf_{u\in \mathcal{N}}I(u)=\hat{m}_0\ge D_{\hat{\sigma}_0}>0. $}
 \par
\vskip2mm
 \noindent
 {\bf Proof.} The proof is similar to Lemma 4.6. Applying (\ref{fc2}) with $t=0$, one has
  \begin{eqnarray}
\label{fc18}
\dfrac{1}{2}f(\upsilon)\upsilon-F(\upsilon)-\dfrac{p_0}{2(\gamma+2)}|\upsilon|^{\gamma+2}\geq0,\;\;\;\; \forall \upsilon\in\mathbb{R}.
 \end{eqnarray}
Following from $(\phi_1)$, $(\phi_2)$, $(\phi_5)$,  $(\ref{b2}),$  $(\ref{b3}),$  $(\ref{ad1})$, $(\ref{c23})$, (\ref{fc18}) and H\"older's inequality, we obtain that for $\forall u\in \mathcal{N}$,
 \begin{eqnarray*}
\label{fc25}
I(u)&=&I(u)-\dfrac{1}{2}\langle I'(u),u\rangle\nonumber\\
&=&\int_{\Omega}\Phi\left(\dfrac{u^2+|\nabla u|^2}{2}\right)dx-\int_{\Omega}F(u)dx-\dfrac{1}{2}\int_{\Omega}\phi\left(\dfrac{u^2+|\nabla u|^2}{2}\right)(u^2+|\nabla u|^2)dx+\dfrac{1}{2}\int_{\Omega}f(u)u dx\nonumber\\
&=&\int_{\Omega}\int^{(u^2+|\nabla u|^2)/2}_{0}\phi(\tau)d\tau dx
-\dfrac{1}{2}\int_{\Omega}\phi\left(\dfrac{u^2+|\nabla u|^2}{2}\right)(u^2+|\nabla u|^2)dx
+ \int_{\Omega}\left[\dfrac{1}{2}f(u)u - F(u)\right]dx\nonumber\\
&=&\int_{\Omega}\phi(\xi_{u,x})\dfrac{u^2+|\nabla u|^2}{2}dx
-\dfrac{1}{2}\int_{\Omega}\phi\left(\dfrac{u^2+|\nabla u|^2}{2}\right)(u^2+|\nabla u|^2)dx
+ \int_{\Omega}\left[\dfrac{1}{2}f(u)u - F(u)\right]dx\nonumber\\
&\geq&\int_{\Omega}\phi\left(\dfrac{u^2+|\nabla u|^2}{2}\right)\dfrac{u^2+|\nabla u|^2}{2}dx
-\dfrac{1}{2}\int_{\Omega}\phi\left(\dfrac{u^2+|\nabla u|^2}{2}\right)(u^2+|\nabla u|^2)dx
+ \int_{\Omega}\left[\dfrac{1}{2}f(u)u - F(u)\right]dx\nonumber\\
&=& \int_{\Omega}\left[\dfrac{1}{2}f(u)u -F(u)\right]dx\nonumber\\
&\geq&\dfrac{p_0\gamma}{2(\gamma+2)}\int_{\Omega}|u|^{\gamma+2}dx\nonumber\\
&\geq&\dfrac{p_0\gamma}{2(\gamma+2)}[\text{vol}(\Omega)]^{1-(\gamma+2)/p}\|u\|_p^{\gamma+2}\nonumber\\
&=&\dfrac{p_0\gamma}{2(\gamma+2)}[\text{vol}(\Omega)]^{1-(\gamma+2)/p}
\left(\dfrac{(1-\sigma_0)\phi_0r^{p-1}}{C_0r^{p-1}+C_0}\right)^{(\gamma+2)/p}
\left(\dfrac{(1-\sigma_0)\phi_0r^{p-1}}{(C_0r^{p-1}+C_0)(C_p)^p}\right)^{[2(\gamma+2)]/[p(p-2)]}\nonumber\\
&=:& \hat{D}_{\sigma_0}> 0.
 \end{eqnarray*}
The proof is completed.
 \qed

\vskip4mm
 \noindent
{\bf Acknowledgement}\\
This work is supported by Yunnan Fundamental Research Projects of China (grant No: 202301AT070465) and  Xingdian Talent
Support Program for Young Talents of Yunnan Province in China.
\vskip2mm
 \noindent
 {\bf Conflict of interest}\\
On behalf of all authors, the  authors states that there is no conflict of interest.\\

\bibliographystyle{amsplain}

 \end{document}